# Global Stability Notions to Enhance the Rigor and Robustness of Adaptive Control

**Iasson Karafyllis[*] and Miroslav Krstic[**]**

[*]Dept. of Mathematics, National Technical University of Athens,
Zografou Campus, 15780, Athens, Greece,
email: iasonkar@central.ntua.gr

[**]Dept. of Mechanical and Aerospace Eng., University of California, San Diego, La Jolla, CA 92093-0411, U.S.A., email: krstic@ucsd.edu

## Abstract

Stability theory plays a crucial role in feedback control. However, adaptive control theory requires advanced and specialized stability notions that are not frequently used in standard feedback control theory. The present document is a set of notes for a graduate course. It describes the global stability notions needed in (robust) adaptive control and develops the mathematical tools that are used for the proof of such stability properties. Moreover, the document shows why and how these global stability properties arise in adaptive control. We focus on stability properties for time-invariant systems. Consequently, tracking control problems are not covered by the present document.

It is a regrettable reality of the adaptive control literature that many authors use terms like "global asymptotic stability" in a loose way and one needs to carefully scrutinize the steps in their analysis in order to understand what the authors actually prove/mean. It is far preferable for adaptive control theory to use a precise terminology that describes the obtained stability properties in a rigorous (and scientifically correct) way. That is why we believe that the present set of notes can inspire and help young researchers in adaptive control.

We try to present the material in the least technical way. Emphasis is placed on the detailed presentation of simple and pedagogical examples. For the interested and mathematically oriented reader we give the precise references where some of the results are proved. We have also included some novel stability results with their proofs. Moreover, we try to present some important issues in adaptive control by means of easy-to-grasp examples.

# Table of Contents

**Notation.** Throughout this document, we adopt the following notation.

∗ $\mathbb{R}_+ := [0, +\infty)$. $x^+ = \max(x, 0)$ for all $x \in \mathbb{R}$.

∗ Let $S \subseteq \mathbb{R}^n$ be an open set and let $A \subseteq \mathbb{R}^n$ be a set that satisfies $S \subseteq A \subseteq cl(S)$, where $cl(S)$ is the closure of $S \subseteq \mathbb{R}^n$. By $C^0(A; \Omega)$, we denote the class of continuous functions on $A$, which take values in $\Omega \subseteq \mathbb{R}^m$. By $C^k(A; \Omega)$, where $k \geq 1$ is an integer, we denote the class of functions on $A \subseteq \mathfrak{R}^n$, which takes values in $\Omega \subseteq \mathbb{R}^m$ and has continuous derivatives of order $k$. In other words, the functions of class $C^k(A; \Omega)$ are the functions which have continuous derivatives of order $k$ in $S = \text{int}(A)$ that can be continued continuously to all points in $\partial S \cap A$. When $\Omega = \mathbb{R}$ then we write $C^0(A)$ or $C^k(A)$.

∗ By $K$ we denote the class of strictly increasing $C^0$ functions $a : \mathbb{R}_+ \to \mathbb{R}_+$ with $a(0) = 0$. By $K_\infty$ we denote the class of strictly increasing $C^0$ functions $a : \mathbb{R}_+ \to \mathbb{R}_+$ with $a(0) = 0$ and $\lim_{s \to +\infty} a(s) = +\infty$. By $KL$ we denote the set of all continuous functions $\sigma : \mathbb{R}_+ \times \mathbb{R}_+ \to \mathbb{R}_+$ with the

properties: (i) for each $t \geq 0$ the mapping $\sigma(\cdot, t)$ is of class $K$; (ii) for each $s \geq 0$, the mapping $\sigma(s, \cdot)$ is non-increasing with $\lim_{s \to +\infty} \sigma(s, t) = 0$.

∗ For a vector $x \in \mathbb{R}^n$, $|x|$ denotes its Euclidean norm and $x'$ denotes its transpose.

∗ Let $S \subseteq \mathbb{R}^n$ be a non-empty set with $0 \in S$. We say that a function $V: S \to \mathbb{R}_+$ is positive definite if $V(x) > 0$ for all $x \in S$ with $x \neq 0$ and $V(0) = 0$. We say that a continuous function $V: S \to \mathbb{R}_+$ is radially unbounded if the following property holds: "for every $M > 0$ the set $\{x \in S : V(x) \leq M\}$ is compact". For $V \in C^1(S; \mathbb{R}_+)$ we define $\nabla V(x) = \left(\frac{\partial V}{\partial x_1}(x), \ldots, \frac{\partial V}{\partial x_n}(x)\right)$.

∗ Let $I \subseteq \mathbb{R}$ be an interval with non-empty interior and let $D \subseteq \mathbb{R}^n$. For every $a \geq 1$, $L^a(I; D)$ denotes the set of equivalence classes of Lebesgue measurable functions $f: I \to D$ for which $\|f\|_a = \left(\int_I |f(t)|^a \, dt\right)^{1/a} < +\infty$. $L^\infty(I; D)$ denotes the set of equivalence classes of Lebesgue measurable functions $f: I \to D$ for which $\|f\|_\infty := \operatorname{ess\,sup}_{t \in I}(|f(t)|) < +\infty$. When $D = \mathbb{R}$, we write $L^a(I)$ or $L^\infty(I)$. For a measurable function $\sup_{t \in I}(|f(t)|)$ always denotes $\operatorname{ess\,sup}_{t \in I}(|f(t)|)$. We also use the space $L^\infty_{loc}(\mathbb{R}_+; D)$. This is the space of equivalence classes of Lebesgue measurable functions $f: \mathbb{R}_+ \to D$ for which their restriction on every bounded interval $I \subseteq \mathbb{R}_+$ with non-empty interior is of class $L^\infty(I; D)$. For example, the function $f(t) = \exp(t)$ does not belong to $L^\infty(\mathbb{R}_+)$ but belongs to $L^\infty_{loc}(\mathbb{R}_+; \mathbb{R})$.

A word of note: if the reader is not familiar with the notion of a Lebesgue measurable function, or with expressions like "a.e.=almost everywhere" and quantities like $\operatorname{ess\,sup}_{t \in I}(|f(t)|)$, the reader can still understand the theory presented in this work. All that is needed is to replace everywhere the term "Lebesgue measurable function" by the term "piecewise continuous function", the expression "a.e." by the expression "with the possible exception of a countable number of points of which there are only finitely many in every bounded interval" and the quantity $\operatorname{ess\,sup}_{t \in I}(|f(t)|)$ by the quantity $\sup_{t \in I}(|f(t)|)$. Of course, what we are suggesting is not equivalent from a mathematical point of view (there are Lebesgue measurable functions that are not piecewise continuous). However, with such a swap of notions the reader can understand the theory without having to learn measure theory.

# 1) Introduction

This set of notes describes the global stability notions needed in (robust) adaptive control and develops the mathematical tools that are used for the proof of such stability properties.

Global stability notions and properties are presented in detail in standard textbooks like [9]. For example, notions like Global Attractivity, Lyapunov stability, Global Asymptotic Stability for input-free systems and the notion of Input-to-State Stability for systems with inputs, are explained in detail in [9]. However, all these notions deal with the evolution of the full state of the closed-loop system, which in the case of adaptive control incorporates both the state of the plant and the state of the parameter estimator (as well as possibly the states of a state observer and of filters for the parameter estimator). In adaptive control, the goal of asymptotic performance is the regulation of the full state (or of an output mapping) of the plant, without having additional requirements (except for boundedness) of the parameter estimator. The state of the plant is thus to be regarded as an output of the overall closed-loop adaptive system. Therefore, all the stability notions given in [9] have to be extended to the case of a system with an output, and such extensions are provided in this document. Another issue that plays important role in adaptive control is the fact that we need to have stability properties that do not show convergence to a point or to a manifold but rather show convergence to a neighborhood of a point or a manifold. Furthermore, it is also necessary to know how large this neighborhood is. This implies the use of practical stability properties.

The task of the development of practical output stability notions has been performed by many researchers in control theory (see notes and comments below). The present set of notes provides and analyzes a collection of stability properties for input-free systems: Lagrange output stability, Lyapunov Output stability, global asymptotic output stability (GAOS), (practical) uniform global asymptotic output stability (p-UGAOS) and (practical) global output attractivity (p-GOA). All these stability notions become the well-known stability notions for input-free systems when the output coincides with the state of the system. Moreover, we also explain in detail a set of stability properties for systems with inputs: (practical) Input-to-Output Stability (p-IOS), the (practical) output asymptotic gain (p-OAG) property, the (practical) Uniform Bounded-Input-Bounded-State (p-UBIBS) property and the zero (practical) Output Asymptotic Gain (zero p-OAG) property.

However, the present set of notes is not simply a set of notes on global output stability properties. We also do the following things:

i) we show how the need of the output stability properties arises in adaptive control theory, and

ii) for each stability notion we provide verifiable sufficient conditions together with simple academic examples that illustrate the use of the provided sufficient conditions for the proof of each stability property.

It is a regrettable reality of the adaptive control literature that many authors use terms like "global asymptotic stability" in a loose way and one needs to carefully scrutinize the steps in their analysis in order to understand what the authors actually prove/mean. For example, there is confusion between global asymptotic stability and global asymptotic output stability. There are papers where the authors prove boundedness of solutions (with no uniformity with respect to the initial conditions) and global output attractivity and still the authors talk about global asymptotic stability. We hope the present set of notes can help the researchers to adopt a precise terminology that describes the obtained stability properties in a rigorous (and scientifically correct) way.

## 2) Global Stability Notions for Input-Free Systems

Let $f : \mathbb{R}^n \to \mathbb{R}^n$ be a locally Lipschitz vector field with $f(0) = 0$ and $h : \mathbb{R}^n \to \mathbb{R}^p$ be a continuous mapping with $h(0) = 0$. Consider the dynamical system

$$\dot{x} = f(x), \, x \in \mathbb{R}^n \tag{1}$$

with output

$$y = h(x) \tag{2}$$

We assume that the dynamical system (1) is forward complete, i.e., for every $x_0 \in \mathbb{R}^n$ the unique solution $x(t) = \phi(t, x_0)$ of the initial-value problem (1) with initial condition $x(0) = x_0$ exists for all $t \geq 0$. We use the notation $y(t, x_0) = h(\phi(t, x_0))$ for all $t \geq 0$, $x_0 \in \mathbb{R}^n$ and $B_R = \{ x \in \mathbb{R}^n : |x| < R \}$ for all $R > 0$. We say that system (1), (2) is

i) *Lagrange output stable* if for every $R > 0$ the set $\{ |y(t, x_0)| : x_0 \in B_R, t \geq 0 \}$ is bounded.

ii) *Lyapunov output stable* if for every $\varepsilon > 0$ there exists $\delta(\varepsilon) > 0$ such that for all $x_0 \in B_{\delta(\varepsilon)}$, it holds that $|y(t, x_0)| \leq \varepsilon$ for all $t \geq 0$.

iii) *Globally Asymptotically Output Stable (GAOS)* if system (1), (2) is Lagrange and Lyapunov output stable and $\lim_{t \to +\infty} (y(t, x_0)) = 0$ for all $x_0 \in \mathbb{R}^n$.

iv) *Uniformly Globally Asymptotically Output Stable (UGAOS)* if system (1), (2) is Lagrange and Lyapunov output stable and for every $\varepsilon, R > 0$ there exists $T(\varepsilon, R) > 0$ such that for all $x_0 \in B_R$, it holds that $|y(t, x_0)| \leq \varepsilon$ for all $t \geq T(\varepsilon, R)$.

We say that system (1), (2) is *Globally Output Attractive (GOA)* if $\lim_{t \to +\infty} (y(t, x_0)) = 0$ for all $x_0 \in \mathbb{R}^n$. We say that system (1), (2) is *practically Globally Output Attractive (p-GOA)* if there exists a constant $\tilde{\alpha} > 0$ such that $\limsup_{t \to +\infty} (|y(t, x_0)|) \leq \tilde{\alpha}$ for all $x_0 \in \mathbb{R}^n$. The constant $\tilde{\alpha} > 0$ is called the *asymptotic residual constant*.

When $h(x) = x$ then the word "output" in the above properties is omitted (e.g., Lagrange stability, Lyapunov stability, GAS, UGAS, p-UGAS, GA, p-GA).

It should be noted that (see Theorem 2.2 on page 62 in [4]) UGAOS for system (1), (2) is equivalent to the existence of a function $\beta \in KL$ such that the following estimate holds for all $x_0 \in \mathbb{R}^n$ and $t \geq 0$:

$$|y(t, x_0)| \leq \beta(|x_0|, t) \tag{3}$$

UGAOS is a much stronger stability property than GAOS. UGAOS guarantees that there are no initial conditions for which convergence of the output to zero happens at a rate that may slow towards zero rate.

We say that system (1), (2) is *practically Uniformly Globally Asymptotically Output Stable (p-UGAOS)* if there exists a function $\beta \in KL$ and a constant $\alpha > 0$ such that the following estimate holds for all $x_0 \in \mathbb{R}^n$ and $t \geq 0$:

$$|y(t, x_0)| \leq \beta(|x_0|, t) + \alpha \tag{4}$$

The constant $\alpha > 0$ is called the *residual constant*. Notice that p-UGAOS with residual constant $\alpha > 0$ is equivalent to UGAOS for system (1) with output $y = (|h(x)| - \alpha)^+$. If system (1), (2) is p-UGAOS with residual constant $\alpha > 0$ then system (1), (2) is p-GOA with asymptotic residual constant $\tilde{\alpha} > 0$ less than or equal to the residual constant $\alpha > 0$, i.e., $\limsup_{t \to +\infty}(|y(t, x_0)|) \leq \tilde{\alpha}$ for certain $\tilde{\alpha} \in (0, \alpha]$.

Lemma 1 in [5] shows that system (1), (2) is Lagrange output stable and Lyapunov output stable if and only if there exists a function $\zeta \in K_\infty$ such that the following estimate holds for all $x_0 \in \mathbb{R}^n$ and $t \geq 0$:

$$|y(t, x_0)| \leq \zeta(|x_0|) \tag{5}$$

None of the above stability properties guarantees or assumes that the set $\{x \in \mathbb{R}^n : h(x) = 0\}$ is invariant for (1). Invariance of the set $\{x \in \mathbb{R}^n : h(x) = 0\}$ for (1) is an additional property that can control the overshoot of the output (but there is no guarantee for this).

## 3) Proving Stability Properties for Input-Free Systems

The most convenient way of proving stability properties is by using Lyapunov-like functions. In this section, we focus on providing sufficient Lyapunov-like conditions for proving the output stability properties that were described above.

A requirement of all stability properties is the property of forward completeness. The following theorem provides the means of proving forward completeness for system (1). It is a special case of Theorem 1 in [7].

**Theorem 1:** *Let $H \in C^1(\mathbb{R}^n; \mathbb{R}_+)$ be a function for which there exist constants $c, \sigma \geq 0$ such that the following inequality holds for all $x \in \mathbb{R}^n$:*

$$\dot{H}(x) = \nabla H(x) f(x) \leq c H(x) + \sigma \tag{6}$$

*Furthermore, for each $r > \inf\{H(x) : x \in \mathbb{R}^n\}$ define the set*

$$S_r = \{x \in \mathbb{R}^n : H(x) \leq r\} \tag{7}$$

*and suppose that $f(S_r)$ is bounded for every $r > \inf\{H(x) : x \in \mathbb{R}^n\}$. Then system (1) is forward complete.*

**Discussion of Theorem 1:** Theorem 1 in the same spirit with the results in [1] for systems defined on $\mathbb{R}^n$ (it is a generalization of certain results in [1]) and with the results in [13] for abstract systems defined on normed linear spaces. However, in contrast with the results in [1], Theorem 1 does not require the set $S_r = \{x \in \mathbb{R}^n : H(x) \leq r\}$ to be bounded for all $r > \inf\{H(x) : x \in \mathbb{R}^n\}$, i.e., Theorem 1 does not require $H$ to be radially unbounded. Instead, Theorem 1 requires boundedness of the set $f(S_r) = \{f(x) : x \in \mathbb{R}^n, H(x) \leq r\}$ for all $r > \inf\{H(x) : x \in \mathbb{R}^n\}$; this condition is automatically satisfied if $S_r = \{x \in \mathbb{R}^n : H(x) \leq r\}$ is bounded for all $r > \inf\{H(x) : x \in \mathbb{R}^n\}$, i.e., when $H$ is radially unbounded, but can also be satisfied in other cases. The following example illustrates this point.

Now we digress to an example that illustrates Theorem 1 but is not inspired by adaptive control.

**Example 1:** Consider the following system

$$\begin{aligned}
\dot{z}_1 &= v_1 \\
\dot{z}_2 &= v_2 \\
\dot{v}_1 &= -g(z_1 - z_2) \\
\dot{v}_2 &= g(z_1 - z_2) \\
x &= (z_1, z_2, v_1, v_2)' \in \mathbb{R}^4
\end{aligned} \tag{8}$$

where $g \in C^1(\mathbb{R})$ is a function that satisfies

$$sg(s) \geq 0, \text{ for all } s \in \mathbb{R} \text{ and } \int_0^{+\infty} g(s)ds = \int_0^{-\infty} g(s)ds = +\infty. \tag{9}$$

Define the function

$$H(x) = \frac{1}{2}v_1^2 + \frac{1}{2}v_2^2 + \int_0^{z_1-z_2} g(s)ds \tag{10}$$

Definition (10) and equations (8) give us for all $x = (z_1, z_2, v_1, v_2)' \in \mathbb{R}^4$:

$$\dot{H}(x) = -v_1 g(z_1 - z_2) + v_2 g(z_1 - z_2) + (v_1 - v_2)g(z_1 - z_2) = 0 \tag{11}$$

Notice that the function $H$ defined by (10) is not radially unbounded. Therefore, we cannot exploit the results in [1] in order to prove forward completeness for system (8). On the other hand, the function $H$ satisfies all requirements of Theorem 1. Indeed, exploiting (9) and (10) we get $\inf\{H(x): x \in \mathbb{R}^4\} = 0$. Moreover, (9) implies that for every $r > 0$ there exists $M(r) > 0$ such that the following implication holds:

$$\int_0^{z_1-z_2} g(s)ds \leq r \Rightarrow |z_1 - z_2| \leq M(r)$$

Thus, since $f(x) = (v_1, v_2, -g(z_1 - z_2), g(z_1 - z_2))'$ (recall (8)), we get for every $r > 0$:

$$S_r \subseteq \left\{(z_1, z_2, v_1, v_2)' \in \mathbb{R}^4 : |v_1| \leq \sqrt{2r}, |v_2| \leq \sqrt{2r}, |z_1 - z_2| \leq M(r)\right\}$$

$$f(S_r) \subseteq \left\{w \in \mathbb{R}^4 : |w| \leq 2\sqrt{2r} + 2\max_{|s| \leq M(r)}(|g(s)|)\right\}$$

where $S_r$ is defined by (7). Since (11) implies that (6) holds with $c = \sigma = 0$, we conclude from Theorem 1 that system (8) is forward complete. ◁

The following proposition provides sufficient Lyapunov-like conditions for Lagrange output stability and Lyapunov output stability of system (1), (2). Its proof is provided in [5].

**Proposition 1 (Lyapunov conditions for Lagrange and Lyapunov output stability):** *Suppose that (1) is forward complete and there exist functions $a \in K_\infty$, $W \in C^1(\mathbb{R}^n; \mathbb{R}_+)$ with $W(0) = 0$, such that the following inequalities hold for all $x \in \mathbb{R}^n$:*

$$a(|h(x)|) \leq W(x) \tag{12}$$

$$\nabla W(x) f(x) \leq 0 \tag{13}$$

*Then system (1), (2) is Lagrange output stable and Lyapunov output stable.*

**Example 2:** Consider again the forward complete system (8) with output

$$y = h(x) = (v_1, v_2)' \in \mathbb{R}^2 \qquad (14)$$

Notice that the function $H$ defined by (10) satisfies the inequality $\frac{1}{2}|h(x)|^2 \leq H(x)$ for all $x \in \mathbb{R}^4$. Therefore, (12) holds with $W = H$ and $a(s) = \frac{1}{2}s^2$. Moreover, (10) implies that $H(0) = 0$ and (11) shows that (13) holds. Therefore, by virtue of Proposition 1 we conclude that system (8), (14) is Lagrange output stable and Lyapunov output stable. ◁

Additional assumptions are needed in order to be able to assert GAOS. The following theorem provides sufficient conditions for GAOS. Its proof is provided in [5].

**Theorem 2 (Lyapunov conditions for GAOS):** *Suppose that (1) is forward complete and there exist functions $a \in K_\infty$, $V, W \in C^1(\mathbb{R}^n; \mathbb{R}_+)$ with $V(0) = W(0) = 0$ and a continuous, positive definite function $\rho: \mathbb{R}_+ \to \mathbb{R}_+$ for which the following inequalities hold for all $x \in \mathbb{R}^n$:*

$$a(|h(x)|) \leq V(x) \qquad (15)$$

$$\dot{V}(x) = \nabla V(x) f(x) \leq -\rho(W(x)) \qquad (16)$$

*Moreover, suppose that (12) holds and that there exists a continuous function $\gamma: \mathbb{R}_+ \to \mathbb{R}_+$ such that one of the following inequalities holds for all $x \in \mathbb{R}^n$:*

$$\dot{W}(x) = \nabla W(x) f(x) \leq \gamma(V(x)) \qquad (17)$$

*or*

$$\dot{W}(x) = \nabla W(x) f(x) \geq -\gamma(V(x)) \qquad (18)$$

*Finally, suppose that either $\rho$ is non-decreasing or that there exists a continuous function $\zeta: \mathbb{R}_+ \to \mathbb{R}_+$ such that the following inequality holds for all $x \in \mathbb{R}^n$:*

$$W(x) \leq \zeta(V(x)) \qquad (19)$$

*Then system (1), (2) is GAOS.*

**Remark:** If $V \in C^1(\mathbb{R}^n; \mathbb{R}_+)$ is a radially unbounded function, then there exist continuous functions $\gamma, \zeta: \mathbb{R}_+ \to \mathbb{R}_+$ such that (17), (18) and (19) hold. Indeed, we obtain in this case the following corollary of Theorem 2, which can also be obtained as a direct consequence of the Barbashin-Krasovskii-LaSalle's theorem.

**Corollary 1:** *Suppose that there exist functions $a \in K_\infty$, $V, W \in C^1(\mathbb{R}^n; \mathbb{R}_+)$ with $V(0) = W(0) = 0$ and a continuous, positive definite function $\rho : \mathbb{R}_+ \to \mathbb{R}_+$ for which (12), (16) hold for all $x \in \mathbb{R}^n$. Moreover, suppose that $V \in C^1(\mathbb{R}^n; \mathbb{R}_+)$ is a radially unbounded function. Then system (1), (2) is GAOS.*

Theorem 2 is less demanding than Barbashin-Krasovskii-LaSalle's theorem. The following example illustrates this point. Both the example and Theorem 2 itself draw some inspiration from adaptive control but go beyond adaptive control.

**Example 3 (Verification of GAOS):** Consider the following three-dimensional system

$$\dot{y} = -(1+w^2)y + \frac{2zg(z,w)}{(1+z^2)^2}$$
$$\dot{z} = -g(z,w)y \tag{20}$$
$$\dot{w} = w + |y|$$
$$x = (y, z, w) \in \mathbb{R}^3$$

where $g : \mathbb{R}^2 \to \mathbb{R}$ is a bounded locally Lipschitz function. We start by showing forward completeness of system (20). This can be shown by using the function $H(x) = (y^2 + z^2 + w^2)/2$ which satisfies the following inequalities for all $x = (y, z, w) \in \mathbb{R}^3$

$$\dot{H}(x) = \nabla H(x) \dot{x} = -(1+w^2)y^2 + \frac{2zg(z,w)}{(1+z^2)^2}y$$
$$- zg(z,w)y + w^2 + w|y| \tag{21}$$
$$\leq -y^2 + 3|g(z,w)||z||y| + w^2 + |w||y|$$
$$\leq 9R^2 z^2 + 3w^2/2 \leq (18R^2 + 3)H(x)$$

where $R := \sup\{|g(z,w)| : (z,w) \in \mathbb{R}^2\}$. It follows from (21) and Theorem 1 (notice that $H$ is radially unbounded) that system (20) is forward complete. Next, we define the function $V(x) = \frac{1}{2}y^2 + \frac{z^2}{1+z^2}$, which satisfies (15) with $a(s) = s^2/2$ for all $s \geq 0$. Moreover, $V$ satisfies the following inequality for all $x = (y, z, w) \in \mathbb{R}^3$

$$\dot{V}(x) = \nabla V(x) \dot{x} = -(1+w^2)y^2 \leq -y^2 \tag{22}$$

Inequality (22) shows that (16) holds with $\rho(s) = 2s$ for all $s \geq 0$ (a non-decreasing function) and $W(x) = y^2/2$. Furthermore, (12) holds with $a(s) = s^2/2$ for all $s \geq 0$. Finally, we show that (17) holds. Indeed, by using the fact that $R := \sup\{|g(z,w)| : (z,w) \in \mathbb{R}^2\}$, we obtain for all $x = (y, z, w) \in \mathbb{R}^3$

$$\dot{W}(x) = \nabla W(x)\dot{x} = -\left(1+w^2\right)y^2 + \frac{2zg(z,w)}{\left(1+z^2\right)^2}y \quad (23)$$

$$\leq -y^2 + R|y| \leq R^2/4$$

Inequality (23) shows that (17) holds with $\gamma(s) \equiv R^2/4$. Notice that since $\rho$ is non-decreasing, we do not have to show (19) (although it happens to hold with $\zeta(s) = s$ for all $s \geq 0$). Therefore Theorem 2 guarantees that system (20) is GAOS.

Barbashin-Krasovskii-LaSalle's theorem could not have been used in order to prove GAOS for (20), because it admits unbounded solutions. To see this notice that the differential inequality $\dot{w} \geq w$ holds for all $x = (y, z, w) \in \mathbb{R}^3$. Thus, we get $w(t) \geq \exp(t)w(0)$ for all $t \geq 0$, and consequently, when $w(0) > 0$ the solutions of (20) are unbounded.

In the same way, Barbălat's Lemma could not have been easily used in order to prove GAOS for (20), because (22) shows that $\dot{W}(x)$ is bounded from above but not necessarily bounded from below. In other words, we cannot (easily) prove that the mapping $t \mapsto W(\phi(t,x))$ is uniformly continuous. ◁

The following example draws inspiration from adaptive control.

**Example 4 (Verification of GAOS):** Consider the system

$$\dot{\xi} = -\left(z + k + \sigma(\theta+z)^2\right)\xi$$
$$\dot{z} = \xi^2 \quad (24)$$
$$x = (\xi, z)' \in \mathbb{R}^2$$

where $\theta \in \mathbb{R}$, $k > 0$, $\sigma \geq 0$ are constants, with output $y = \xi$. Using the radially unbounded function

$$V(x) = \frac{1}{2}\xi^2 + \frac{1}{2}z^2 \quad (25)$$

that satisfies for all $x \in \mathbb{R}^2$ the equation

$$\dot{V}(x) = -\left(k + \sigma(\theta+z)^2\right)\xi^2 \leq -k\xi^2 \leq 0 \quad (26)$$

we are in a position to conclude that:

a) by virtue of Theorem 1, system (24) is forward complete,
b) the set $\{(0, z)' \in \mathbb{R}^2, z \in \mathbb{R}\}$ is invariant for (24); it is the set of equilibrium points of system (24),
c) by virtue of Proposition 1, system (24) is Lagrange stable and Lyapunov stable,
d) by using Corollary 1 with $W(x) = \frac{1}{2}\xi^2$ and $\rho(s) = 2ks$, system (24) is GAOS.

However, we notice that the convergence rate of $y = \xi \in \mathbb{R}$ to zero is not necessarily uniform with respect to the initial condition, i.e., UGAOS is not necessarily valid for (24). To see this, consider the case $\sigma = 0$ for which we are able to write a formula for the solution of (24) with initial condition $(\xi(0), z(0)) = (\xi_0, z_0) \in \mathbb{R}^2$:

$$\xi(t) = \frac{K\xi_0}{K\cosh(Kt) + (z_0 + k)\sinh(Kt)}$$

$$z(t) = -k - K + \frac{2K(z_0 + k + K)\exp(2Kt)}{(z_0 + k + K)(\exp(2Kt) - 1) + 2K}$$

where $K := \sqrt{\xi_0^2 + (k + z_0)^2}$. It is clearly shown by the above formulas that the convergence rate is state-dependent. Notice that convergence can happen at a rate that may slow towards zero rate. Indeed, when $z_0 = -k$ then $\xi(t) = \frac{2\xi_0 \exp(-|\xi_0|t)}{1 + \exp(-2|\xi_0|t)}$, which shows that the convergence rate tends to zero as $\xi_0$ tends to zero. Moreover, following the arguments in [17], we are able to show that system (24) with $\sigma = 0$ is only GAOS, and not UGAOS. We repeat the argument here: when $\xi_0 > 0$ and $z_0 + k < 0$, it holds that $\xi(t^*) = \max_{t \geq 0}(\xi(t)) = \sqrt{\xi_0^2 + (k + z_0)^2}$ for

$$t^* = \frac{1}{2\sqrt{\xi_0^2 + (k + z_0)^2}} \ln\left(\frac{\sqrt{\xi_0^2 + (k + z_0)^2} - z_0 - k}{\sqrt{\xi_0^2 + (k + z_0)^2} + z_0 + k}\right) > 0.$$ If system (24) with $\sigma = 0$ were UGAOS, then there would exist $T > 0$ such that for all $x_0 = (\xi_0, z_0) \in \mathbb{R}^2$ with $|x_0| \leq \sqrt{1 + (k+1)^2}$ it holds that $|y(t, x_0)| \leq \frac{1}{2}$ for all $t \geq T$. Clearly, the sequence of initial conditions $x_{0,i} = \left(\frac{1}{i}, -k - 1\right) \in \mathbb{R}^2$ with $i \geq 1$ being an integer, satisfies $|x_{0,i}| \leq \sqrt{1 + (k+1)^2}$ and $y(t_i^*, x_{0,i}) = \xi(t_i^*) = \sqrt{1 + i^{-2}} > \frac{1}{2}$ for

$$t_i^* = \frac{1}{2\sqrt{1 + i^{-2}}} \ln\left(1 + 2i^2\left(\sqrt{1 + i^{-2}} + 1\right)\right) \geq \frac{\ln(1 + 4i^2)}{2\sqrt{2}}.$$ Consequently, we would get $T \geq t_i^* \geq \frac{\ln(1 + 4i^2)}{2\sqrt{2}}$ for all integers $i \geq 1$; a contradiction. Thus, system (24) with $\sigma = 0$ is not UGAOS.

The example is continued below. ◁

The counterexample in [17] showed that the existence of functions $V \in C^1(\mathbb{R}^n; \mathbb{R}_+)$ with $V(0) = 0$, $a, c \in K_\infty$ for which (15) holds and for which the following inequality holds for all $x \in \mathbb{R}^n$

$$\dot{V}(x) = \nabla V(x) f(x) \leq -c(|h(x)|) \tag{27}$$

cannot guarantee UGAOS. On the other hand, it is known that the existence of functions $V \in C^1(\mathbb{R}^n; \mathbb{R}_+)$ with $V(0) = 0$, $a \in K_\infty$ and a continuous, positive definite function $\rho : \mathbb{R}_+ \to \mathbb{R}_+$ for which (15) holds and for which the following inequality holds for all $x \in \mathbb{R}^n$

$$\dot{V}(x) = \nabla V(x) f(x) \leq -\rho(V(x)) \tag{28}$$

does guarantee UGAOS without any further assumption (see for instance the discussion in [23, 24] as well as Theorem 2.4 on page 90 in [4]). The key difference between (27) and (28) is that the dissipation rate in the former involves only the output norm, whereas the latter involves the Lyapunov function itself. However, in practice it is difficult to obtain a dissipation rate, like in (28), that depends on the Lyapunov function itself.

The following theorem provides Lyapunov-like conditions for UGAOS of system (1), (2) without requiring a dissipation rate that involves the Lyapunov function and shows the missing link between GAOS and UGAOS in terms of the dissipation rate. Its proof is provided in [5].

**Theorem 3 (Lyapunov conditions for UGAOS):** *Suppose that (1) is forward complete and there exist functions $V, W \in C^1(\mathbb{R}^n; \mathbb{R}_+)$, $a \in K_\infty$ and a continuous, positive definite function $\rho : \mathbb{R}_+ \to \mathbb{R}_+$ such that inequalities (12), (13), (16) hold (which are given again below for reader's convenience)*

$$a(|h(x)|) \leq W(x) \tag{12}$$

$$\dot{W}(x) = \nabla W(x) f(x) \leq 0 \tag{13}$$

$$\dot{V}(x) = \nabla V(x) f(x) \leq -\rho(W(x)) \tag{16}$$

*Then system (1), (2) is UGAOS.*

The above statement shows that a sufficient additional requirement in order for (16) and (12) to ensure UAOS is that the dissipation rate $W$ does not increase along the system's solutions.

**Remark:** Theorem 3 does not assume that (15) holds. Moreover, it is not assumed in Theorem 3 that $V(0) = 0$. Notice that if (15) and (28) hold then all assumptions of Theorem 3 are valid with $W = V$, in which case we recover a classical sufficient condition for UGAOS; see [4, 23, 24].

**Example 4 (continued-verification of UGAOS):** Consider again system (24) with $\sigma, k > 0$. We next show that system (24) with $k \geq \frac{1}{4\sigma} + \theta$ and output $y = \xi$ is UGAOS. Using the radially unbounded function defined by (25) that satisfies for all $x \in \mathbb{R}^2$ equation (26), we notice that (26) implies (16) with $W(x) = \xi^2 / 2$ and $\rho(s) = 2ks$. It is also clear that (12) holds with $a(s) = s^2 / 2$. Finally, we obtain for all $x \in \mathbb{R}^2$:

$$\dot{W}(x) = -\left(z + k + \sigma(\theta + z)^2\right)\xi^2 \tag{29}$$

Since the inequality $z + k + \sigma(\theta + z)^2 \geq k - \frac{1}{4\sigma} - \theta$ holds for all $z \in \mathbb{R}$, we conclude from (29) that (13) holds when $k \geq \frac{1}{4\sigma} + \theta$. Therefore, Theorem 3 implies that system (24) is UGAOS when $k \geq \frac{1}{4\sigma} + \theta$.

We also notice that equation (29) implies the differential inequality $\dot{W}(x) \leq -2\left(k - \frac{1}{4\sigma} - \theta\right)W(x)$ for all $x \in \mathbb{R}^2$. Integrating the above differential inequality, we obtain the estimate

$$|y(t, x_0)| = |\xi(t)| \leq \exp\left(-\left(k - \frac{1}{4\sigma} - \theta\right)t\right)|\xi_0|$$

that holds for all $x_0 = (\xi_0, z_0)' \in \mathbb{R}^2$, $t \geq 0$. Therefore, when $k > \frac{1}{4\sigma} + \theta$ we are able to conclude that the output converges to zero at an exponential rate (notice that uniform attractivity does not necessarily imply an exponential convergence rate). The above estimate is independent of $z_0$ and that is the reason that such a strong stability property is termed in [24] as "State-Independent Uniform Output Stability". However, the reader should notice that a state-independent stability estimate (i.e., an estimate independent of $z_0$) cannot be easily shown in the case $k = \frac{1}{4\sigma} + \theta$, where we know that UGAOS holds.

Estimates like the ones provided above can be combined with other estimates in order to prove stability properties. Indeed, using (26) we can derive the estimate $\xi^2(t) + z^2(t) \leq \xi^2(0) + z^2(0)$ for all $t \geq 0$. Combining this estimate with the estimate $|\xi(t)| \leq \exp\left(-\left(k - \frac{1}{4\sigma} - \theta\right)t\right)|\xi(0)|$ that holds when $k > \frac{1}{4\sigma} + \theta$, we can guarantee the following estimate for all $t \geq 0$ when $k > \frac{1}{4\sigma} + \theta$:

$$|\xi(t)z(t)| \leq \exp\left(-\left(k - \frac{1}{4\sigma} - \theta\right)t\right)|\xi(0)||x(0)|$$
$$\leq \exp\left(-\left(k - \frac{1}{4\sigma} - \theta\right)t\right)|x(0)|^2$$

The above estimate is an estimate of the form (3) and shows UGAOS for system (24) with $k > \frac{1}{4\sigma} + \theta$ and output $y = \xi z$. ◁

## 4) The Need of Output Stability Properties for Disturbance-Free Systems in Adaptive Control

Adaptive control is a particular case in which output stability notions are useful. Adaptive control strategies rely on a state extension to provide a dynamical estimate of unknown parameters. In this context, it is not always necessary to precisely estimate these unknown parameters to get a satisfactory behavior of the system's solutions. In other words, only part of the state variables of the closed-loop system are requested to converge to zero, which can be rephrased as a GAOS objective. To show this, consider the simplest of cases: the so-called case of a system that satisfies the "matching condition". Consider the following control system

$$\dot{y} = f(y) + g(y)\left(u + (\varphi(y))' \theta\right)$$
$$y \in \mathbb{R}^n, \theta \in \mathbb{R}^p, u \in \mathbb{R} \tag{30}$$

where $f, g : \mathbb{R}^n \to \mathbb{R}^n$, $\varphi : \mathbb{R}^n \to \mathbb{R}^p$ are locally Lipschitz vector fields with $f(0) = 0$ and $\varphi(0) = 0$, $y \in \mathbb{R}^n$ denotes the state, $u \in \mathbb{R}$ is the control input and $\theta \in \mathbb{R}^p$ is the constant vector of unknown parameters. The control system (30) satisfies the "matching condition" (following the terminology in [10]) since the uncertain parameters are in the span of the control.

The following assumption is used for the control system (30).

**(H)** *There exist functions $P \in C^2(\mathbb{R}^n; \mathbb{R}_+)$ being positive definite and radially unbounded, $Q \in C^1(\mathbb{R}^n; \mathbb{R}_+)$ being positive definite and a locally Lipschitz function $k \in C^0(\mathbb{R}^n; \mathbb{R})$ with $k(0) = 0$ such that*

$$\nabla P(y)f(y) + \nabla P(y)g(y)k(y) \le -Q(y) \text{, for all } y \in \mathbb{R}^n \tag{31}$$

When Assumption (H) holds then it is possible to regulate the state to zero by applying a standard adaptive control scheme. To see this, consider the adaptive controller

$$u = k(y) - (\varphi(y))' \hat{\theta}$$
$$\frac{d\hat{\theta}}{dt} = \Gamma^{-1} \nabla P(y) g(y) \varphi(y) \tag{32}$$

where $\Gamma \in \mathbb{R}^{p \times p}$ is a constant, positive definite matrix (parameters of the controller). System (30) in closed loop with (32) then reads:

$$\dot{y} = f(y) + g(y)\left(k(y) - (\varphi(y))' z\right)$$
$$\dot{z} = \Gamma^{-1} \nabla P(y) g(y) \varphi(y) \tag{33}$$
$$x = (y, z) \in \mathbb{R}^n \times \mathbb{R}^p$$

with $z = \hat{\theta} - \theta$ and output map defined by
$$h(x) = y \tag{34}$$

Using the radially unbounded function $V(x) := P(y) + \frac{1}{2}z'\Gamma z$, we obtain (by virtue of (31), (33)) the following differential inequality for all $x = (y, z) \in \mathbb{R}^n \times \mathbb{R}^p$:

$$\dot{V}(x) \leq -Q(y) \quad (35)$$

Using (35) we are in a position to conclude that:

a) by virtue of Theorem 1, system (33) is forward complete,
b) the set $\{(0, z)' \in \mathbb{R}^n \times \mathbb{R}^p, z \in \mathbb{R}^p\}$ is invariant for (33); it is a set of equilibrium points for system (33),
c) by virtue of Proposition 1, system (33) is Lagrange stable and Lyapunov stable,
d) by using Corollary 1 with $W(x) = P(y)$, system (33) is GAOS.

The reader should notice the precision that the stability properties of system (33) are described by using the output stability notions that we introduced above. The above set of stability properties provides important performance guarantees for the closed-loop system; this is the reason for the success of adaptive control.

However, the convergence rate of $y \in \mathbb{R}^n$ to zero is not necessarily uniform with respect to the initial condition, i.e., UGAOS is not necessarily valid for (33). This is a source of various problems that can arise for the closed-loop system (33). The following sections will shed light to the problems that can appear for (33).

## 5) Global Stability Notions for Systems with (Disturbance) Inputs

Let $D \subseteq \mathbb{R}^p$ be a given closed set with $0 \in D$ and $f : \mathbb{R}^n \times D \to \mathbb{R}^n$ be a locally Lipschitz with respect to $x \in \mathbb{R}^n$ mapping with $f(0,0) = 0$. Consider the control system

$$\dot{x} = f(x,d), \; x \in \mathbb{R}^n, \; d \in D \tag{36}$$

We assume that system (36) is forward complete, i.e., for every $x_0 \in \mathbb{R}^n$ and for every Lebesgue measurable and locally essentially bounded input $d : \mathbb{R}_+ \to D$ the unique solution $x(t) = \phi(t, x_0; d)$ of the initial-value problem (36) with initial condition $x(0) = x_0$ corresponding to input $d : \mathbb{R}_+ \to D$ exists for all $t \geq 0$. We use the notation $y(t, x_0; d) = h(\phi(t, x_0; d))$ for all $t \geq 0$, $x_0 \in \mathbb{R}^n$ and for every Lebesgue measurable and locally essentially bounded input $d : \mathbb{R}_+ \to D$.

We say that system (36), (2) is *Input-to-Output Stable (IOS)* if there exist a function $\beta \in KL$ and a non-decreasing, continuous function $\gamma : \mathbb{R}_+ \to \mathbb{R}_+$ with $\gamma(0) = 0$ such that the following estimate holds for all $x_0 \in \mathbb{R}^n$, $t \geq 0$ and for every $d \in L^\infty(\mathbb{R}_+; D)$:

$$|y(t, x_0; d)| \leq \beta(|x_0|, t) + \gamma(\|d\|_\infty) \tag{37}$$

We say that system (36), (2) is *practically Input-to-Output Stable (p-IOS)* if there exist a function $\beta \in KL$, a non-decreasing, continuous function $\gamma : \mathbb{R}_+ \to \mathbb{R}_+$ with $\gamma(0) = 0$ and a constant $\alpha > 0$ such that the following estimate holds for all $x_0 \in \mathbb{R}^n$, $t \geq 0$ and for every $d \in L^\infty(\mathbb{R}_+; D)$:

$$|y(t, x_0; d)| \leq \beta(|x_0|, t) + \gamma(\|d\|_\infty) + \alpha \tag{38}$$

The constant $\alpha > 0$ is called the *residual constant* while the function $\gamma$ is called the *gain function of the input* $d \in D$ *to the output* $y$. It is clear that the gain function is a measure of the sensitivity of system (36), (2) to the disturbance $d \in D$. Notice that p-IOS with residual constant $\alpha > 0$ is equivalent to IOS for system (36) with output $y = (|h(x)| - \alpha)^+$. In the disturbance-free case, the p-IOS property coincides with the p-UGAOS property.

We say that system (36), (2) satisfies the *practical Output Asymptotic Gain (p-OAG)* property if there exists a non-decreasing, continuous function $\tilde{\gamma} : \mathbb{R}_+ \to \mathbb{R}_+$ with $\tilde{\gamma}(0) = 0$ and a constant $\tilde{\alpha} > 0$ such that the following estimate holds for all $x_0 \in \mathbb{R}^n$ and for every $d \in L^\infty(\mathbb{R}_+; D)$:

$$\limsup_{t \to +\infty} (|y(t, x_0; d)|) \leq \tilde{\gamma}(\|d\|_\infty) + \tilde{\alpha} \tag{39}$$

The constant $\tilde{\alpha} > 0$ is called the *asymptotic residual constant* while the non-decreasing, continuous function $\tilde{\gamma}$ is called the *asymptotic gain function of the input* $d \in D$ *to the output* $y$. It is clear that the asymptotic gain function and the asymptotic residual constant are additional measures of the sensitivity of system (36), (2) to the disturbance $d \in D$. It is a consequence of the semigroup property that if system (36), (2) satisfies the p-IOS property then system (36), (2) also satisfies the p-OAG property with asymptotic gain function being less than or equal to the gain function (i.e., $\tilde{\gamma}(s) \leq \gamma(s)$ for all $s \geq 0$) and asymptotic residual constant being less than or equal to the residual constant (i.e., $\tilde{\alpha} \leq \alpha$). In the disturbance-free case, the p-OAG property coincides with the p-GOA property.

When $\tilde{\gamma} \equiv 0$ we say that system (36), (2) satisfies the *zero practical Output Asymptotic Gain property (zero p-OAG)*.

When $h(x) = x$ then the word "output" in the above properties is either replaced by the word "state" (e.g., ISS, p-ISS) or is omitted (e.g., p-AG, zero p-AG).

We say that system (36) satisfies the *practical Uniform Bounded-Input-Bounded-State (p-UBIBS)* property if there exists a function $\bar{\gamma} \in K_\infty$ and a constant $\bar{\alpha} > 0$ such that the following estimate holds for all $x_0 \in \mathbb{R}^n$ and for every $d \in L^\infty(\mathbb{R}_+; D)$:

$$\sup_{t \geq 0}\left(|\phi(t, x_0; d)|\right) \leq \bar{\gamma}\left(|x_0|\right) + \bar{\gamma}\left(\|d\|_\infty\right) + \bar{\alpha} \tag{40}$$

The p-UBIBS property with $\bar{\alpha} = 0$ is called the UBIBS property.

Clearly, the p-UBIBS property is equivalent to the existence of a continuous function $B: \mathbb{R}^n \times \mathbb{R}_+ \to \mathbb{R}_+$ for which the following estimate holds for all $x_0 \in \mathbb{R}^n$ and for every $d \in L^\infty(\mathbb{R}_+; D)$:

$$\sup_{t \geq 0}\left(|\phi(t, x_0; d)|\right) \leq B\left(x_0, \|d\|_\infty\right) \tag{41}$$

In the disturbance-free case, the p-UBIBS property coincides with the Lagrange stability property. Moreover, in the disturbance-free case, the UBIBS property coincides with the combination of Lagrange and Lyapunov stability.

## 6) Proving Stability Properties for Systems with (Disturbance) Inputs

A requirement of all stability properties for systems with inputs is the property of forward completeness. The following theorem is an extension of Theorem 1 and provides the means of proving forward completeness for system (36).

**Theorem 4:** *Let $H \in C^1(\mathbb{R}^n; \mathbb{R}_+)$ be a function for which there exist non-decreasing functions $c, \sigma \in C^0(\mathbb{R}_+; \mathbb{R}_+)$ such that the following inequality holds for all $(x, d) \in \mathbb{R}^n \times D$:*

$$\dot{H}(x,d) = \nabla H(x) f(x,d) \leq c(|d|) H(x) + \sigma(|d|) \tag{42}$$

*Furthermore, for each $r > \inf\{H(x) : x \in \mathbb{R}^n\}$ define by means of (7) the set $S_r$ as well as the set*

$$D_r = \{d \in D : |d| \leq r\} \tag{43}$$

*and suppose that $f(S_r \times D_r)$ is bounded for every $r > \inf\{H(x) : x \in \mathbb{R}^n\}$. Then system (36) is forward complete.*

**Proof:** Let $x_0 \in \mathbb{R}^n$ and $d \in L^\infty_{loc}(\mathbb{R}_+; D)$ be given. Since $f$ is locally Lipschitz with respect to $x \in \mathbb{R}^n$, the initial-value problem $\dot{x} = f(x, d)$ with initial condition $x(0) = x_0$ has a unique solution $x(t) \in \mathbb{R}^n$ defined for $t \in [0, t_{\max})$, where $t_{\max} \in (0, +\infty]$. Furthermore, if $t_{\max} < +\infty$ then $\limsup_{t \to t_{\max}^-}(|x(t)|) = +\infty$. We show next by means of a contradiction argument that $t_{\max} = +\infty$.

Without loss of generality, we may assume that (42) holds with $\sigma \equiv 0$ (if not then we can replace $H$ by $1 + H$ and $c$ by $c + \sigma$). Suppose that $t_{\max} < +\infty$ and define $\|d\| = \sup_{0 \leq t \leq t_{\max}}(|d(t)|)$.

Inequality (42) and definition $\|d\| = \sup_{0 \leq t \leq t_{\max}}(|d(t)|)$ as well as the fact that $c$ is non-decreasing, give the following differential inequality for $t \in [0, t_{\max})$ a.e.:

$$\frac{d}{dt}(H(x(t))) \leq c(\|d\|) H(x(t)) \tag{44}$$

Integrating (44), we obtain that

$$H(x(t)) \leq \exp(c(\|d\|) t) H(x_0), \text{ for all } t \in [0, t_{\max}) \tag{45}$$

which implies that

$$H(x(t)) \leq \exp(c(\|d\|) t_{\max}) H(x_0), \text{ for all } t \in [0, t_{\max}) \tag{46}$$

Let $r = \exp(c(\|d\|)t_{\max})H(x_0) + \|d\|$ and notice that definitions $\|d\| = \sup_{0 \leq t \leq t_{\max}}(|d(t)|)$, (7), (43) and estimate (46) imply that $(x(t), d(t)) \in S_r \times D_r$ for $t \in [0, t_{\max})$ a.e.. Since by assumption, $f(S_r \times D_r)$ is bounded, it follows that there exists a constant $M > 0$, such that

$$|\dot{x}(t)| \leq M, \text{ for } t \in [0, t_{\max}) \text{ a.e.} \tag{47}$$

Inequality (47) implies that

$$|x(t) - x_0| \leq M t, \text{ for all } t \in [0, t_{\max})$$

and consequently (using the triangle inequality)

$$|x(t)| \leq |x_0| + M t_{\max}, \text{ for all } t \in [0, t_{\max}) \tag{48}$$

Therefore, estimate (48) shows that $\limsup_{t \to t_{\max}^-}(|x(t)|) = +\infty$ cannot hold. Consequently, $t_{\max} = +\infty$. The proof is complete. ◁

**Discussion of Theorem 4:** Theorem 4 in the same spirit with the results in [1] for systems defined on $\mathbb{R}^n$ (it is a generalization of certain results in [1]) and with the results in [13] for abstract systems defined on normed linear spaces. However, in contrast with the results in [1], Theorem 4 does not require the set $S_r = \{x \in \mathbb{R}^n : H(x) \leq r\}$ to be bounded for all $r > \inf\{H(x) : x \in \mathbb{R}^n\}$, i.e., Theorem 4 does not require $H$ to be radially unbounded. Instead, Theorem 4 requires boundedness of the set $f(S_r \times D_r)$ for all $r > \inf\{H(x) : x \in \mathbb{R}^n\}$; this condition is automatically satisfied if $H$ is radially unbounded (notice that in this case $S_r = \{x \in \mathbb{R}^n : H(x) \leq r\}$ is compact for all $r > \inf\{H(x) : x \in \mathbb{R}^n\}$ and since $D \subseteq \mathbb{R}^p$ is closed it also follows that $D_r = \{d \in D : |d| \leq r\}$ is compact for all $r > \inf\{H(x) : x \in \mathbb{R}^n\}$).

**Example 5:** Consider the following system

$$\begin{aligned} \dot{x} &= f(x, d) \\ x &\in \mathbb{R}^n, d \in \mathbb{R}^p \end{aligned} \tag{49}$$

where $f \in C^1(\mathbb{R}^n \times \mathbb{R}^p)$ is a mapping that satisfies the following estimate for all $x \in \mathbb{R}^n, d \in \mathbb{R}^p$

$$|f(x, d)| \leq (L + \kappa(|d|))|x| + \kappa(|d|) \tag{50}$$

where $\kappa \in K_\infty$ and $L > 0$ is a constant. Define the radially unbounded function

$$H(x) = \frac{1}{2}|x|^2 \tag{51}$$

for which we get (using (50) and completing the squares) for all $x \in \mathbb{R}^n, d \in \mathbb{R}^p$

$$\begin{aligned}\dot{H}(x,d) &= \nabla H(x)f(x,d) = x'f(x,d) \\ &\leq |x||f(x,d)| \leq \left(L+\kappa(|d|)\right)|x|^2 + |x|\kappa(|d|) \\ &\leq \left(L+\kappa(|d|)+\frac{1}{2}\right)|x|^2 + \frac{1}{2}\kappa^2(|d|) \\ &= \left(2L+2\kappa(|d|)+1\right)H(x) + \frac{1}{2}\kappa^2(|d|)\end{aligned} \quad (52)$$

It follows from (52) that (42) holds for appropriate non-decreasing functions $c, \sigma \in C^0(\mathbb{R}_+; \mathbb{R}_+)$. We conclude from Theorem 4 that system (49) is forward complete. ◁

The p-UBIBS property is a more demanding property than simple forward completeness because the p-UBIBS property requires boundedness of solutions for all bounded inputs. Therefore, there are systems which are forward complete but do not satisfy the p-UBIBS property (e.g., the system $\dot{x} = d$ with $x \in \mathbb{R}, d \in \mathbb{R}$). The following result provides sufficient conditions for the p-UBIBS property.

**Theorem 5 (Sufficient Lyapunov-like conditions for the p-UBIBS property):** *Let $H_i \in C^1(\mathbb{R}^n; \mathbb{R}_+)$, $i = 1, ..., m$ be a family of functions for which there exists a non-decreasing function $R \in C^0(\mathbb{R}_+; \mathbb{R}_+)$ such that the following implication holds for all $(x,d) \in \mathbb{R}^n \times D$ and $i = 1, ..., m$:*

$$H_i(x) \geq R(|d|) \Rightarrow \dot{H}_i(x,d) = \nabla H_i(x)f(x,d) \leq 0 \quad (53)$$

*Furthermore, suppose that the function $H(x) = \max_{i=1,...,m}(H_i(x))$ is a radially unbounded function. Then system (36) satisfies the p-UBIBS property. Moreover, if in addition the function $H(x) = \max_{i=1,...,m}(H_i(x))$ is positive definite and $R(0) = 0$ then system (36) satisfies the UBIBS property.*

**Proof:** Let arbitrary $x_0 \in \mathbb{R}^n$ and $d \in L^\infty(\mathbb{R}_+; D)$ be given. Since $f$ is locally Lipschitz with respect to $x \in \mathbb{R}^n$, the initial-value problem $\dot{x} = f(x,d)$ with initial condition $x(0) = x_0$ has a unique solution $x(t) \in \mathbb{R}^n$ defined for $t \in [0, t_{\max})$, where $t_{\max} \in (0, +\infty]$. Furthermore, if $t_{\max} < +\infty$ then $\limsup_{t \to t_{\max}^-}(|x(t)|) = +\infty$.

We claim that the following estimate holds for all $t \in [0, t_{\max})$ and $i = 1, ..., m$:

$$H_i(x(t)) \leq \max\left(H_i(x_0), R(\|d\|_\infty)\right) \quad (54)$$

The proof is made by contradiction. Suppose that there exists $i=1,...,m$ and $T \in [0, t_{max})$ such that $H_i(x(T)) > \max(H_i(x_0), R(\|d\|_\infty))$. Since $H_i(x(0)) \leq \max(H_i(x_0), R(\|d\|_\infty))$, by continuity of the mapping $t \mapsto H_i(x(t))$ the set $\{t \in [0,T]: H_i(x(t)) = \max(H_i(x_0), R(\|d\|_\infty))\}$ is non-empty. Define $\tau = \sup\{t \in [0,T]: H_i(x(t)) = \max(H_i(x_0), R(\|d\|_\infty))\}$ and notice that continuity of the mapping $t \mapsto H_i(x(t))$ implies that $H_i(x(\tau)) = \max(H_i(x_0), R(\|d\|_\infty))$ and $\tau < T$. Definition of $\tau$ guarantees that $H_i(x(t)) \geq \max(H_i(x_0), R(\|d\|_\infty))$ for all $t \in [\tau, T]$. Thus, we obtain from (53) and the fact that $R \in C^0(\mathbb{R}_+; \mathbb{R}_+)$ is a non-decreasing function for $t \in [\tau, T]$ a.e.:

$$\frac{d}{dt}(H_i(x(t))) \leq 0 \qquad (55)$$

Inequality (55) implies that $H_i(x(T)) \leq H_i(x(\tau)) = \max(H_i(x_0), R(\|d\|_\infty))$, which contradicts the estimate $H_i(x(T)) > \max(H_i(x_0), R(\|d\|_\infty))$.

The fact that $H(x) = \max_{i=1,...,m}(H_i(x))$ and (54) implies that the following estimate holds for all $t \in [0, t_{max})$:

$$H(x(t)) \leq \max(H(x_0), R(\|d\|_\infty)) \qquad (56)$$

Since $H \in C^0(\mathbb{R}^n; \mathbb{R}_+)$ is a radially unbounded function we can guarantee from (56) that $t_{max} = +\infty$.

Since $H \in C^0(\mathbb{R}^n; \mathbb{R}_+)$ is a radially unbounded function the minimum $b = \min\{H(x): x \in \mathbb{R}^n\}$ is well-defined and satisfies the inequality $b \geq 0$. We define the non-decreasing function $a(s) := \max\{|x|: H(x) \leq s\}$ for $s \geq b$, which is well-defined by virtue of the fact that $H$ is a radially unbounded function. Then, we get $|x| \leq a(H(x))$ for all $x \in \mathbb{R}^n$. The function $\tilde{a}(s) := \int_s^{s+1} a(p)dp$ is a continuous, non-decreasing function defined on $[b, +\infty)$ that satisfies $a(s) \leq \tilde{a}(s)$ for all $s \geq b$. Consequently, we get $|x| \leq \tilde{a}(H(x))$ for all $x \in \mathbb{R}^n$. Using (56) we get the following estimate for all $t \geq 0$:

$$|x(t)| \leq \tilde{a}(\max(H(x_0), R(\|d\|_\infty)))$$

It follows from the above estimate that (41) holds with $B(x,s) = \tilde{a}(\max(H(x), R(s)))$ for all $(x,s) \in \mathbb{R}^n \times \mathbb{R}_+$. Therefore, the p-UBIBS property holds for (36).

If in addition the function $H(x) = \max_{i=1,\ldots,m}(H_i(x))$ is positive definite and $R(0) = 0$ then Lemma 3.5 on page 138 in [9] implies the existence of functions $a_1, a_2 \in K_\infty$ such that the following inequality holds for all $x \in \mathbb{R}^n$:

$$a_1(|x|) \leq H(x) \leq a_2(|x|)$$

Combining the above inequality with (56) we obtain the following estimate for all $t \geq 0$:

$$|x(t)| \leq \max\left(a_1^{-1}(a_2(|x_0|)), a_1^{-1}(R(\|d\|_\infty))\right)$$
$$\leq a_1^{-1}(a_2(|x_0|)) + a_1^{-1}(R(\|d\|_\infty))$$

The above estimate shows the validity of (40) with $\bar{\alpha} = 0$ and $\bar{\gamma}(s) = a_1^{-1}(a_2(s)) + a_1^{-1}(R(s))$. Since $\bar{\gamma} \in K_\infty$ (recall that $R(0) = 0$), we conclude that the UBIBS property holds for (36). The proof is complete. ◁

**Example 6:** Consider the following system

$$\begin{aligned}\dot{x}_1 &= p(x,d)\left((1+d-x_1^2)x_1 + d\right) \\ \dot{x}_2 &= -q(x,d)(x_2+d) \\ x &= (x_1, x_2) \in \mathbb{R}^2, d \in \mathbb{R}\end{aligned} \quad (57)$$

where $p, q : \mathbb{R}^2 \times \mathbb{R} \to \mathbb{R}_+$ are non-negative, locally Lipschitz functions. Define the functions

$$H_i(x) = \frac{1}{2}x_i^2, \quad i = 1, 2 \quad (58)$$

for which we have for all $(x,d) \in \mathbb{R}^2 \times \mathbb{R}$:

$$\begin{aligned}\dot{H}_1(x,d) &= \nabla H_1(x)\dot{x} = p(x,d)\left((1+d-x_1^2)x_1^2 + x_1 d\right) \\ &\leq p(x,d)\left((1+|d|)x_1^2 - x_1^4 + |d||x_1|\right) \\ &\leq p(x,d)\left(2(1+|d|)H_1(x) - 4H_1^2(x) + \frac{1}{2}|d|^2 + H_1(x)\right) \\ &\leq p(x,d)\left(2(|d|+2-2H_1(x))H_1(x) + \frac{1}{2}|d|^2 - H_1(x)\right)\end{aligned} \quad (59)$$

It follows that when $H_1(x) \geq \frac{1}{2}|d|^2$ and $H_1(x) \geq \frac{|d|+2}{2}$ then $\nabla H_1(x)\dot{x} \leq 0$. Moreover, we also have for all $(x,d) \in \mathbb{R}^2 \times \mathbb{R}$:

$$\dot{H}_2(x,d) = \nabla H_2(x)\dot{x} = -q(x,d)\left(x_2^2 + x_2 d\right) \leq -\frac{q(x,d)}{2}\left(2H_2(x) - d^2\right)$$

It follows that when $H_2(x) \geq |d|^2/2$ then $\nabla H_2(x)\dot{x} \leq 0$.

Therefore, (53) holds with $R(s) = \max(s+2, s^2)/2$. Since $H(x) = \max_{i=1,2}(H_i(x)) = \max(x_1^2, x_2^2)/2$ is a radially unbounded function, Theorem 5 implies the p-UBIBS property for system (57). ◁

The following theorem provides sufficient conditions for the p-IOS property. It is a direct consequence of Lemma 2.14 on page 82 and Remark 2.4 on page 85 in [4].

**Theorem 6 (Lyapunov conditions for p-IOS):** *Suppose that (36) is forward complete and there exist a function $a \in K_\infty$, a non-decreasing, continuous function $\varphi: \mathbb{R}_+ \to \mathbb{R}_+$, functions $V_i \in C^1(\mathbb{R}^n; \mathbb{R}_+)$ $i = 1,...,m$ with $V_i(0) = 0$ for $i = 1,...,m$ and continuous, positive definite functions $\rho_i: \mathbb{R}_+ \to \mathbb{R}_+$, $i = 1,...,m$ for which the following inequality holds for all $x \in \mathbb{R}^n$*

$$a(|h(x)|) \leq \max_{i=1,...,m}(V_i(x)) \tag{60}$$

*and the following implication holds for all $x \in \mathbb{R}^n$, $d \in D$, $i = 1,...,m$:*

$$V_i(x) \geq \varphi(|d|) \Rightarrow \dot{V}_i(x,d) = \nabla V_i(x) f(x,d) \leq -\rho_i(V_i(x)) \tag{61}$$

*Then system (36), (2) satisfies the p-IOS property with gain function $\gamma(s) = a^{-1}(\varphi(s)) - a^{-1}(\varphi(0))$ and residual constant $\alpha = a^{-1}(\varphi(0))$. Moreover, if $\varphi(0) = 0$ then system (36), (2) satisfies the IOS property with gain function $\gamma = a^{-1} \circ \varphi$.*

**Example 7:** Consider the following system

$$\begin{aligned}\dot{y} &= -\left(z + k + \sigma(\theta + z)^2\right)y + d \\ \dot{z} &= y^2 \\ x &= (y,z)' \in \mathbb{R}^2, d \in \mathbb{R}\end{aligned} \tag{62}$$

where $\theta \in \mathbb{R}$, $\sigma, k > 0$ are constants that satisfy the inequality $k > \theta + 1/(4\sigma)$. The disturbance-free version of system (62) was studied in Example 4. The radially unbounded function

$$H(x) = y^2/2 + z^2/2 \tag{63}$$

satisfies for all $x \in \mathbb{R}^2, d \in \mathbb{R}$ the equation

$$\dot{H}(x,d) = \nabla H(x) f(x,d) = -\left(k + \sigma(\theta + z)^2\right) y^2 + yd \leq -\frac{k}{2}y^2 + \frac{1}{2k}d^2 \tag{64}$$

Inequality (64) implies that all assumptions of Theorem 4 hold for system (62). Therefore, system (62) is forward complete. The function

$$V(x) = y^2/2 \tag{65}$$

satisfies for all $x \in \mathbb{R}^2, d \in \mathbb{R}$ the inequality (using the fact that $z + k + \sigma(\theta + z)^2 \geq k - \theta - 1/(4\sigma)$ for all $z \in \mathbb{R}$):

$$\dot{V}(x,d) = \nabla V(x) f(x,d) = -\left(z + k + \sigma(\theta + z)^2\right) y^2 + dy$$
$$\leq -\left(k - \frac{1}{4\sigma} - \theta\right) y^2 + |d||y| \tag{66}$$

The function $V$ defined by (65) satisfies inequality (60) with $a(s) = s^2/2$. Moreover, due to (66) we obtain the implication

$$\frac{1}{2\varepsilon^2}|d|^2 \leq V(x) \Rightarrow \dot{V}(x,d) = \nabla V(x) f(x,d) \leq -\varepsilon y^2 = -2\varepsilon V(x) \tag{67}$$

where

$$\varepsilon := \frac{1}{2}\left(k - \frac{1}{4\sigma} - \theta\right) > 0 \tag{68}$$

Consequently, implication (61) holds with $\varphi(s) = s^2/(2\varepsilon^2)$, $\rho(s) = 2\varepsilon s$ and $r = 0$. We conclude from Theorem 6 that (62) satisfies the IOS property with gain function $\gamma(s) = s/\varepsilon^2$. ◁

Small-Gain arguments can also be used for the derivation of IOS estimates (see for instance [4]). Here we prefer to give a (vector) Lyapunov-like version of the small-gain theorem instead of the usual trajectory-based version.

**Theorem 7 (Small-Gain conditions for IOS):** *Suppose that there exist a function $a \in K_\infty$, non-decreasing, continuous functions $\gamma_1, \gamma_2, \zeta, p : \mathbb{R}_+ \to \mathbb{R}_+$ with $\gamma_1(0) = \gamma_2(0) = \zeta(0) = p(0) = 0$, functions $V_1, V_2, U \in C^1(\mathbb{R}^n; \mathbb{R}_+)$ with $V_1(0) = V_2(0) = 0$ and $U$ being positive definite and radially unbounded, a constant $c \geq 0$ and a continuous, positive definite function $\rho : \mathbb{R}_+ \to \mathbb{R}_+$ for which the following inequalities hold for all $x \in \mathbb{R}^n$, $d \in D$*

$$a(|h(x)|) \leq \max(V_1(x), V_2(x)) \tag{69}$$

$$\dot{U}(x,d) = \nabla U(x) f(x,d) \leq cU(x) + p(V_1(x)) + p(V_2(x)) + \zeta(|d|) \tag{70}$$

*and the following implications hold for all $x \in \mathbb{R}^n$, $d \in D$:*

$$V_1(x) \geq \max\left(\zeta(|d|), \gamma_1(V_2(x))\right) \Rightarrow \dot{V}_1(x,d) = \nabla V_1(x) f(x,d) \leq -\rho(V_1(x)) \tag{71}$$

$$V_2(x) \geq \max\left(\zeta(|d|), \gamma_2(V_1(x))\right) \Rightarrow \dot{V}_2(x,d) = \nabla V_2(x) f(x,d) \leq -\rho(V_2(x)) \tag{72}$$

*Moreover, suppose that $\gamma_1(\gamma_2(s)) < s$ and $\gamma_2(\gamma_1(s)) < s$ for all $s > 0$. Then system (36), (2) satisfies the IOS property.*

**Proof:** The proof is a simple application of Theorem 5.3 on page 219 in [4]. More, specifically we apply Theorem 5.3 in [4] with $W(t,x) = \exp(-(c+1)t)U(x)$. Due to the fact that $U$ is positive definite and radially unbounded, there exist functions $a_1, a_2 \in K_\infty$ such that $a_1(|x|) \leq U(x) \leq a_2(|x|)$ for all $x \in \mathbb{R}^n$. Using Lemma 3.2 in [4], we obtain a function $q \in K_\infty$ such that $a_1^{-1}(rs) \leq q(r)q(s)$ for all $r, s \geq 0$. Defining $\mu(t) = \dfrac{1}{q(\exp((c+1)t))}$ we conclude that the inequality $q^{-1}(\mu(t)s) \leq \exp(-(c+1)t)a_1(s)$ holds for all $t, s \geq 0$. Putting everything together, we get the inequality $q^{-1}(\mu(t)|x|) \leq W(t,x) \leq a_2(|x|)$ for all $t \geq 0$, $x \in \mathbb{R}^n$. Using (70) we also obtain the inequality for all $t \geq 0$, $x \in \mathbb{R}^n$ and $d \in D$:

$$\frac{\partial W}{\partial t}(t,x) + \frac{\partial W}{\partial x}(t,x)f(x,d) \leq -W(t,x) + p(V_1(x)) + p(V_2(x)) + \zeta(|d|)$$

The proof is complete. ◁

**Example 8:** Consider the following system

$$\begin{aligned} \dot{x}_1 &= -x_1 + bx_2^m + d \\ \dot{x}_2 &= -x_2^m + x_1 \\ x &= (x_1, x_2)' \in \mathbb{R}^2, d \in \mathbb{R} \end{aligned} \quad (73)$$

where $b \in \mathbb{R}$ is a constant with $|b| < 1$ and $m$ is an odd positive integer. Since $|b| < 1$ there exists $\lambda \in (0,1)$ with $|b| < \lambda^2$. The functions

$$V_i(x) = x_i^2 / 2, \; i = 1, 2 \quad (74)$$

satisfy inequality (69) with $h(x) = x$ and $a(s) = s^2 / 4$. Moreover, we have the following implications for all $x \in \mathbb{R}^2, d \in \mathbb{R}$:

$$\max\left(\frac{2}{(1-\lambda)^2}|d|^2, \frac{2^{m-1}b^2}{\lambda^2}V_2^m(x)\right) \leq V_1(x) \Rightarrow \dot{V}_1(x,d) = \nabla V_1(x)\dot{x} \leq -(1-\lambda)V_1(x) \quad (75)$$

$$\frac{2^{(1-m)/m}}{\lambda^{2/m}} V_1^{1/m}(x) \leq V_2(x) \Rightarrow \dot{V}_2(x,d) = \nabla V_2(x)\dot{x} \leq -2(1-\lambda)V_2(x) \quad (76)$$

It follows that implications (71), (72) hold with $\rho(s) = (1-\lambda)s$, $\gamma_1(s) = 2^{m-1}b^2 s^m / \lambda^2$, $\gamma_2(s) = 2^{(1-m)/m} s^{1/m} / \lambda^{2/m}$ and $\zeta(s) = 2s^2 / (1-\lambda)^2$ for $s \geq 0$. The small-gain conditions $\gamma_1(\gamma_2(s)) < s$ and $\gamma_2(\gamma_1(s)) < s$ for all $s > 0$ hold when $|b| < \lambda^2$. Finally, the positive definite and radially unbounded function $U(x) = V_1(x) + V_2(x)$ satisfies for all $x \in \mathbb{R}^2, d \in \mathbb{R}$:

$$\dot{U}(x,d) = \nabla U(x)f(x,d) = -x_1^2 + bx_1 x_2^m + dx_1 - x_2^{m+1} + x_1 x_2$$

$$\leq \frac{b^2 + 1}{2} x_1^2 + \frac{1}{2} x_2^{2m} + \frac{1}{2} x_2^2 + \frac{2}{(1-\lambda)^2}|d|^2 + \frac{(1-\lambda)^2}{8} x_1^2$$

The above inequality shows that (70) holds for appropriate $p \in K_\infty$ and appropriate $c \geq 0$. We conclude from Theorem 7 that system (73) is ISS. ◁

The zero p-OAG property is a novel property that has not been studied in the literature. The only result that allows the proof of the zero p-OAG property deals with a specific class of systems: systems with deadzone. The following theorem is the foundation upon which we build the robust adaptive designs in [6, 8]. It employs what at first might seem like an unusual form of a deadzone acting on the Lyapunov function, with an update gain $\exp(-z(t))$ that decays as the gain $z(t)$ being updated increases. With the analysis tool presented in this theorem we are able to construct robust adaptive designs which achieve the following properties: practical regulation to an arbitrarily small set for the plant state, regulation to an even smaller set if the asymptotic values of the disturbance and the unknown parameter are small, with the set size proportional to the small asymptotic value of the disturbance and parameter, an exponential transient decay of the plant state to a residual set proportional to the "historic" size of the disturbance and parameter, and the prevention of the unbounded drift of the adaptive gain $z(t)$.

**Theorem 8 (Lyapunov conditions for systems with deadzone):** *Let $\Theta \subseteq \mathbb{R}^p$, $D \subseteq \mathbb{R}^l$ be closed sets with $0 \in \Theta$ and $0 \in D$. Consider the system*

$$\dot{\xi} = f(\xi, z, \theta, d)$$
$$\dot{z} = \Gamma \exp(-z)\big(V(\xi, z) - \varepsilon\big)^+ \tag{77}$$

*where $x = (\xi, z) \in \mathbb{R}^n \times \mathbb{R}$ is the state, $\theta \in \Theta \subseteq \mathbb{R}^p$, $d \in D \subseteq \mathbb{R}^l$ are disturbances, $\Gamma, \varepsilon > 0$ are constants, $f : \mathbb{R}^n \times \mathbb{R} \times \Theta \times D \to \mathbb{R}^n$ is a locally Lipschitz with respect to $(x, z) \in \mathbb{R}^n \times \mathbb{R}$ mapping with $f(0, z, \theta, 0) = 0$ for all $(z, \theta) \in \mathbb{R} \times \Theta$. Let $V \in C^1(\mathbb{R}^n \times \mathbb{R})$ with $V(0, z) = 0$ for all $z \in \mathbb{R}$ and $V(\xi, z) > 0$ for all $\xi \neq 0$, $z \in \mathbb{R}$ be a function that satisfies the following properties:*

**(i)** *For every $M > 0$ the set $\{(\xi, z) \in \mathbb{R}^n \times [-M, M] : V(\xi, z) \leq M\}$ is bounded.*

**(ii)** *There exist functions $\kappa, \lambda, \rho \in K_\infty$ and constants $a, c, \delta > 0$, $b \geq 0$, for which the inequality $\rho(s) \geq cs$ holds for all $s \in [0, \delta]$ and such that the following inequality holds for all $(\xi, z) \in \mathbb{R}^n \times \mathbb{R}, \theta \in \Theta, d \in D$:*

$$\frac{\partial V}{\partial \xi}(\xi, z) f(\xi, z, \theta, d) + \Gamma \exp(-z) \frac{\partial V}{\partial z}(\xi, z)\big(V(\xi, z) - \varepsilon\big)^+$$
$$\leq -\rho(V(\xi, z)) + \chi(|d|, |\theta|, \exp(z)) \tag{78}$$

*where $\chi(s_1, s_2, s_3) = a \dfrac{s_1^2 + \big((s_2 - b - \lambda(s_3))^+\big)^2}{1 + \kappa(s_3)}$ for all $s_1, s_2, s_3 \geq 0$.*

*Then there exist functions $B \in C^0(\mathbb{R}^n \times \mathbb{R} \times \mathbb{R}_+ \times \mathbb{R}_+)$, $\gamma \in K_\infty$, $\sigma \in KL$ such that for every $(\xi_0, z_0) \in \mathbb{R}^n \times \mathbb{R}$, $d \in L^\infty(\mathbb{R}_+; D)$, $\theta \in L^\infty(\mathbb{R}_+; \Theta)$ the unique solution of (77) with $(\xi(0), z(0)) = (\xi_0, z_0)$ exists and satisfies the following estimates for all $t \geq 0$:*

$$\limsup_{t \to +\infty}(V(\xi(t), z(t))) \leq \varepsilon, \tag{79}$$

$$\limsup_{t \to +\infty}(V(\xi(t), z(t))) \leq \gamma\left(\chi\left(\limsup_{t \to +\infty}(|d(t)|), \limsup_{t \to +\infty}(|\theta(t)|), \exp\left(\lim_{t \to +\infty}(z(t))\right)\right)\right) \tag{80}$$

$$V(\xi(t), z(t)) \leq \sigma(V(\xi_0, z_0), t) + \gamma\left(\chi\left(\|d\|_\infty, \|\theta\|_\infty, \exp(z_0)\right)\right), \tag{81}$$

$$|\xi(t)| \leq B\left(\xi_0, z_0, \|d\|_\infty, \|\theta\|_\infty\right), \tag{82}$$

$$z_0 \leq z(t) \leq \lim_{s \to +\infty}(z(s)) \leq B\left(\xi_0, z_0, \|d\|_\infty, \|\theta\|_\infty\right). \tag{83}$$

*Moreover, if the inequality* $\rho(s) \geq cs$, *holds for all* $s \geq 0$ *then* $\gamma(s) = c^{-1}s$ *for all* $s \geq 0$ *and* $\sigma(s,t) = \exp(-ct)s$ *for all* $t, s \geq 0$.

**Remark:** Estimate (80) shows that $\lim_{t \to +\infty}(|\xi(t)|) = 0$ whenever $\lim_{t \to +\infty}(d(t)) = 0$ and $\limsup_{t \to +\infty}(|\theta(t)|) \leq b + \lambda\left(\exp\left(\lim_{t \to +\infty}(z(t))\right)\right)$. Inequality (79) shows the zero p-OAG property for system (77) with output $y = V(\xi, z)$ with asymptotic residual constant $\varepsilon > 0$. Inequalities (82), (83) show that system (77) satisfies the p-UBIBS property.

**Proof:** We start by noticing that there exist functions $\sigma \in KL$, $\gamma \in K_\infty$ with the following property: for every $T \in (0, +\infty]$, $\mu \geq 0$ and for every absolutely continuous function $y : [0,T) \to \mathbb{R}_+$ for which the following inequality holds for $t \in [0,T)$ a.e.:

$$\dot{y}(t) \leq -\rho(y(t)) + \mu \tag{84}$$

the following estimate holds for all $t \in [0,T)$:

$$y(t) \leq \sigma(y(0), t) + \gamma(\mu) \tag{85}$$

The existence of functions $\sigma \in KL$, $\gamma \in K_\infty$ follows from Lemma 2.14 on page 82 in [4] and by noticing that (91) shows that the implication

$$y(t) \geq \rho^{-1}(2\mu) \Rightarrow \dot{y}(t) \leq -\rho(y(t))/2$$

holds for $t \in [0,T)$ a.e.. Moreover, if the inequality $\rho(s) \geq cs$, holds for all $s \geq 0$ then by directly integrating (84) we get (85) with $\gamma(s) = c^{-1}s$ for all $s \geq 0$ and $\sigma(s,t) = \exp(-ct)s$ for all $t, s \geq 0$.

Define $\bar{c}(r) := \inf\left\{\frac{\rho(l)}{l} : 0 < l \leq r\right\}$ for all $r > 0$. Notice that continuity of $\rho \in K_\infty$ and the fact that the inequality $\rho(l) \geq cl$ holds for all $l \in [0, \delta]$, guarantee that $\bar{c} : (0, +\infty) \to (0, +\infty)$ is a well-defined, non-increasing function that satisfies

$$\rho(l) \geq \bar{c}(r)l, \text{ for all } r > 0 \text{ and } l \in [0, r] \tag{86}$$

Without loss of generality, we may assume that $\bar{c} : (0, +\infty) \to (0, +\infty)$ is continuous (otherwise replace the function $\bar{c}$ by the function $\tilde{c}(r) := \int_r^{r+1} \bar{c}(l)dl$ defined for $r > 0$). Next define for $r \geq 0$ the non-increasing, continuous and positive function

$$C(r) := \bar{c}\left(1 + \sigma(r, 0) + \gamma(ar)\right) \tag{87}$$

Moreover, define for $r \geq 0$

$$A(r) := \kappa^{-1}\left(\left(\frac{ar}{C(r)\varepsilon} - 1\right)^+\right), \tag{88}$$

$$G(r) := \ln\left(1 + \max(A(r), \exp(r)) + \frac{\Gamma}{C(r)}\left(1 + \frac{a}{C(r)}\right)r\right). \tag{89}$$

Notice that both $A$ and $G$ are non-decreasing, continuous functions.

Let arbitrary $(\xi_0, z_0) \in \mathbb{R}^n \times \mathbb{R}$, $d \in L^\infty(\mathbb{R}_+; D)$, $\theta \in L^\infty(\mathbb{R}_+; \Theta)$ be given. Since the functions $f : \mathbb{R}^n \times \mathbb{R} \times \Theta \times D \to \mathbb{R}^n$ and $V$ are locally Lipschitz with respect to $(x, z) \in \mathbb{R}^n \times \mathbb{R}$ there exists an unique solution $x(t) = (\xi(t), z(t))$ of (77) with $x(0) = (\xi(0), z(0)) = (\xi_0, z_0)$ defined on a maximal interval $[0, t_{max})$, where $t_{max} \in (0, +\infty]$. For every $t \in [0, t_{max})$, we have from (77) that $\dot{z}(t) \geq 0$. Therefore, $z(t)$ is non-decreasing on $[0, t_{max})$. Moreover, the mapping $z \mapsto \chi(|d|, |\theta|, \exp(z))$ is non-increasing. Using (78) we conclude that the following inequalities hold almost everywhere in $[0, t_{max})$:

$$\frac{d}{dt}V(\xi(t), z(t)) \leq -\rho(V(\xi(t), z(t))) + \chi(|d(t)|, |\theta(t)|, \exp(z(t))) \tag{90}$$

which implies

$$\frac{d}{dt}V(\xi(t), z(t)) \leq -\rho(V(\xi(t), z(t))) + \chi(\|d\|_\infty, \|\theta\|_\infty, \exp(z_0)) \tag{91}$$

It follows from (84) and (85) that the following estimate holds for all $t \in [0, t_{max})$:

$$V(\xi(t), z(t)) \leq \sigma(V(\xi_0, z_0), t) + \gamma\left(\chi(\|d\|_\infty, \|\theta\|_\infty, \exp(z_0))\right) \tag{92}$$

Moreover, if the inequality $\rho(s) \geq cs$ holds for all $s \geq 0$ then $\gamma(s) = c^{-1}s$ for all $s \geq 0$ and $\sigma(s, t) = \exp(-ct)s$ for all $t, s \geq 0$.

It follows from (90), (86), estimate (92) and definition (87) that the following inequality holds for $t \in [0, t_{max})$ a.e.:

$$\frac{d}{dt}V(\xi(t), z(t)) \leq -C(s)V(\xi(t), z(t)) + a\frac{\|d\|_\infty^2 + \left((\|\theta\|_\infty - b - \lambda(\exp(z_0)))^+\right)^2}{1 + \kappa(\exp(z(t)))} \tag{93}$$

where
$$s := V(\xi_0, z_0) + |z_0| + \|\theta\|_\infty^2 + \|d\|_\infty^2 \tag{94}$$

Integrating (93) and using the fact that $z(t)$ is non-decreasing on $[0, t_{max})$, we obtain the following estimate for all $t \in [0, t_{max})$:

$$V(\xi(t), z(t)) \leq \exp(-C(s)t) V(\xi_0, z_0) + \frac{1}{C(s)} \chi\left(\|d\|_\infty, \|\theta\|_\infty, \exp(z_0)\right) \tag{95}$$

Moreover, from estimate (95) and definition (94), we get the following estimate for all $t \in [0, t_{max})$:

$$V(\xi(t), z(t)) \leq V(\xi_0, z_0) + \frac{a}{C(s)}\left(\|\theta\|_\infty^2 + \|d\|_\infty^2\right) \leq s + \frac{a}{C(s)} s \tag{96}$$

Next, we show that $z(t)$ is bounded. Since $z(t)$ is non-decreasing on $[0, t_{max})$, we have $z(t) \geq z_0$ for all $t \in [0, t_{max})$. We next distinguish the following cases:

<u>*Case 1:*</u> $\exp(z(t)) \leq 1 + \max(A(s), \exp(s))$ for all $t \in [0, t_{max})$

In this case, it follows from (89) that:
$$z(t) \leq G(s), \text{ for all } t \in [0, t_{max})$$

<u>*Case 2:*</u> There exists $T \in (0, t_{max})$ such that $\exp(z(T)) > 1 + \max(A(s), \exp(s))$. Since $\exp(z(0)) \leq \exp(s) < 1 + \max(A(s), \exp(s))$, continuity of $z(t)$ guarantees that there exist $t_0 \in (0, T)$ such that:
$$\exp(z(t_0)) = 1 + \max(A(s), \exp(s)) \tag{97}$$

Since $z(t)$ is non-decreasing, it follows that:
$$\exp(z(t)) \geq 1 + \max(A(s), \exp(s)), \text{ for } t \in [t_0, t_{max}) \tag{98}$$

Inequality (98) in conjunction with (94) and (88) gives the following inequalities for $t \in [t_0, t_{max})$:

$$\exp(z(t)) \geq A(s) \geq \kappa^{-1}\left(\left(\frac{as}{C(s)\varepsilon} - 1\right)^+\right) \geq \kappa^{-1}\left(\left(a\frac{\|\theta\|_\infty^2 + \|d\|_\infty^2}{C(s)\varepsilon} - 1\right)^+\right)$$

$$\Rightarrow \kappa(\exp(z(t))) \geq \left(a\frac{\|\theta\|_\infty^2 + \|d\|_\infty^2}{C(s)\varepsilon} - 1\right)^+ \geq a\frac{\|\theta\|_\infty^2 + \|d\|_\infty^2}{C(s)\varepsilon} - 1 \Rightarrow C(s)\varepsilon \geq a\frac{\|\theta\|_\infty^2 + \|d\|_\infty^2}{1 + \kappa(\exp(z(t)))}$$

Consequently, we obtain from (93) for almost all $t \in [t_0, t_{max})$:

$$\frac{d}{dt} V(\xi(t), z(t)) \leq -C(s) V(\xi(t), z(t)) + C(s)\varepsilon$$

Integrating the above differential inequality, we obtain for all $t \in [t_0, t_{max})$:

$$V(\xi(t), z(t)) \leq \exp\left(-C(s)(t-t_0)\right)V(\xi(t_0), z(t_0)) + \varepsilon \tag{99}$$

Using (77) and inequality (99), we have for all $t \in [t_0, t_{max})$:

$$\frac{d}{dt}\left(\exp(z(t))\right) \leq \Gamma \exp\left(-C(s)(t-t_0)\right)V(\xi(t_0), z(t_0)) \tag{100}$$

Integrating (100) and using (96) and (97), we obtain for all $t \in [t_0, t_{max})$:

$$\begin{aligned}\exp(z(t)) &\leq \exp\left(z(t_0)\right) + \frac{\Gamma}{C(s)}V(\xi(t_0), z(t_0)) \\ &= 1 + \max\left(A(s), \exp(s)\right) + \frac{\Gamma}{C(s)}V(\xi(t_0), z(t_0)) \\ &\leq 1 + \max\left(A(s), \exp(s)\right) + \frac{\Gamma}{C(s)}\left(1 + \frac{a}{C(s)}\right)s\end{aligned} \tag{101}$$

Since $z(t)$ is non-decreasing, inequality (101) holds for all $t \in [0, t_{max})$.

It follows from (101) and definitions (94), (89) that in both cases the following estimates hold for all $t \in [0, t_{max})$:

$$z_0 \leq z(t) \leq G(s) \tag{102}$$

$$|z(t)| \leq G(s) \tag{103}$$

Estimate (103) implies that $z(t)$ is bounded. Therefore, assumption (i) in conjunction with (92) guarantees that $\xi(t)$ is bounded for all $t \in [0, t_{max})$. Consequently, $t_{max} = +\infty$. Moreover, estimate (92) implies directly the required estimate (81).

Since $z(t)$ is non-decreasing and bounded from above, $\lim_{\tau \to +\infty}\left(z(\tau)\right)$ exists. Thus, we obtain from (102):

$$z_0 \leq z(t) \leq \lim_{\tau \to +\infty}\left(z(\tau)\right) \leq G(s) \tag{104}$$

Next, define the following function for $r \geq 0$:

$$H(r) := \max\left\{|w| : V(w, \zeta) \leq \left(1 + \frac{a}{C(r)}\right)r \text{ and } |\zeta| \leq G(r)\right\} \tag{105}$$

Assumption (i) guarantees that $H$ is well-defined. Since $G$ is non-decreasing and $C$ is non-increasing, it follows from (105) that $H$ is non-decreasing. Moreover, it follows from estimates (96), (103) and definitions (94), (105):

$$|x(t)| \leq H(s), \text{ for all } t \geq 0 \tag{106}$$

Define the continuous function $\tilde{H}(r) := G(r) + \int_{r}^{r+1} H(l)dl$ for $r \geq 0$ (notice that since $H$ is non-decreasing, it is locally bounded and Riemann-integrable on every bounded interval of $\mathbb{R}_+$). We conclude from (106), (102) and (94) that estimates (83), (82) hold with $B(\xi_0, z_0, \|d\|_\infty, \|\theta\|_\infty) = \tilde{H}\left(V(\xi_0, z_0) + |z_0| + \|\theta\|_\infty^2 + \|d\|_\infty^2\right)$.

We next show estimate (79). Estimate (102) and the fact that $z(t)$ is non-decreasing guarantees that the function $z(t)$ has a finite limit as $t \to +\infty$. This implies that the function $\exp(z(t))$ has a finite limit as $t \to +\infty$. Moreover, the facts that $d \in L^\infty(\mathbb{R}_+; D)$, $\theta \in L^\infty(\mathbb{R}_+; \Theta)$ and (103), (106) imply that $\frac{d}{dt}(V(\xi(t), z(t)))$ is of class $L^\infty(\mathbb{R}_+)$. It follows that the function $\frac{d}{dt}(\exp(z(t))) = \Gamma(V(\xi(t), z(t)) - \varepsilon)^+$ is uniformly continuous, in addition to $\int_0^{+\infty} \frac{d}{dt}(\exp(z(t)))dt = \lim_{t \to +\infty}(\exp(z(t))) - \exp(z_0) \leq \exp(G(s)) - \exp(z_0) < +\infty$. From Barbălat's Lemma (see [5, 9]), we have:

$$\lim_{t \to +\infty}\left(\frac{d}{dt}(\exp(z(t)))\right) = \lim_{t \to +\infty}\left(\Gamma(V(\xi(t), z(t)) - \varepsilon)^+\right) = 0 \quad (107)$$

Therefore, estimate (79) holds.

Next, we show estimate (80). Let arbitrary $\epsilon > 0$ be given. Then there exists $\bar{T} > 0$ such that

$$|d(t)| \leq \limsup_{\tau \to +\infty}(|d(\tau)|) + \epsilon \text{ for } t \geq \bar{T} \text{ a.e.}$$

$$|\theta(t)| \leq \limsup_{\tau \to +\infty}(|\theta(\tau)|) + \epsilon \text{ for } t \geq \bar{T} \text{ a.e.}$$

$$z(t) \geq \lim_{\tau \to +\infty}(z(\tau)) - \epsilon \text{ for all } t \geq \bar{T}$$

Define $L_d = \limsup_{\tau \to +\infty}(|d(\tau)|)$, $L_\theta = \limsup_{\tau \to +\infty}(|\theta(\tau)|)$ and $L_z = \lim_{\tau \to +\infty}(z(\tau))$. It follows from (90) that the following differential inequality holds for $t \geq \bar{T}$ a.e.:

$$\frac{d}{dt}V(\xi(t), z(t)) \leq -\rho(V(\xi(t), z(t))) + \chi(L_d + \epsilon, L_\theta + \epsilon, \exp(L_z - \epsilon)) \quad (108)$$

We get from (84), (85) and (108) for all $t \geq \bar{T}$:

$$V(\xi(t), z(t)) \leq \sigma\left(V(\xi(\bar{T}), z(\bar{T})), t - \bar{T}\right) + \gamma\left(\chi(L_d + \epsilon, L_\theta + \epsilon, \exp(L_z - \epsilon))\right) \quad (109)$$

Consequently, we get from estimate (109):

$$\limsup_{t \to +\infty}(V(\xi(t), z(t))) \leq \gamma\left(\chi(L_d + \epsilon, L_\theta + \epsilon, \exp(L_z - \epsilon))\right) \quad (110)$$

Since (110) holds for arbitrary $\epsilon > 0$, we obtain the required estimate (80).
The proof is complete. ◁

**Example 9:** This example conveys the essence of the adaptive robust stability studies in [6, 8]. It covers the closed-loop dynamics for a scalar linear plant. Consider the following system

$$\dot{\xi} = \theta\xi - \left(c + \frac{(1+\exp(z))^2}{2c} + \frac{1+\exp(z)}{4a}(1+\xi^2)\right)\xi + d$$

$$\dot{z} = \Gamma\exp(-z)\left(\frac{1}{2}\xi^2 - \varepsilon\right)^+ \tag{111}$$

$$x = (\xi, z)' \in \mathbb{R}^2, d \in \mathbb{R}, \theta \in \mathbb{R}$$

where $a, c, \Gamma, \varepsilon > 0$ are constants. The function $V(\xi, z) = \xi^2 / 2$ satisfies for all $(\xi, z) \in \mathbb{R}^2$, $\theta, d \in \mathbb{R}$:

$$\frac{\partial V}{\partial \xi}(\xi, z)\dot{\xi} + \frac{\partial V}{\partial z}(\xi, z)\dot{z} = \theta\xi^2 - \left(c + \frac{(1+\exp(z))^2}{2c} + \frac{1+\exp(z)}{4a}(1+\xi^2)\right)\xi^2 + \xi d$$

$$\leq |\theta|\xi^2 - \left(c + \frac{(1+\exp(z))^2}{2c} + \frac{1+\exp(z)}{4a}(1+\xi^2)\right)\xi^2 + \frac{a|d|^2}{1+\exp(z)} + \frac{1+\exp(z)}{4a}\xi^2$$

$$= (|\theta|-1-\exp(z))\xi^2 + (1+\exp(z))\xi^2 - \left(c + \frac{(1+\exp(z))^2}{2c} + \frac{1+\exp(z)}{4a}\xi^2\right)\xi^2 + \frac{a|d|^2}{1+\exp(z)}$$

$$\leq (|\theta|-1-\exp(z))^+ \xi^2 + (1+\exp(z))\xi^2 - \left(c + \frac{(1+\exp(z))^2}{2c} + \frac{1+\exp(z)}{4a}\xi^2\right)\xi^2 + \frac{a|d|^2}{1+\exp(z)}$$

$$\leq a\frac{\left((|\theta|-1-\exp(z))^+\right)^2}{1+\exp(z)} + \frac{1+\exp(z)}{4a}\xi^4 + \frac{c}{2}\xi^2 + \frac{(1+\exp(z))^2}{2c}\xi^2$$

$$-\left(c + \frac{(1+\exp(z))^2}{2c} + \frac{1+\exp(z)}{4a}\xi^2\right)\xi^2 + \frac{a|d|^2}{1+\exp(z)}$$

$$= -\frac{c}{2}\xi^2 + a\frac{|d|^2 + \left((|\theta|-1-\exp(z))^+\right)^2}{1+\exp(z)} = -cV(\xi, z) + a\frac{|d|^2 + \left((|\theta|-1-\exp(z))^+\right)^2}{1+\exp(z)}$$

Consequently, inequality (78) holds with $b = 1$, $\kappa(s) = \lambda(s) = s$ for $s \geq 0$. It follows that system (111) satisfies all assumptions of Theorem 8. Therefore, there exists a function $B \in C^0(\mathbb{R} \times \mathbb{R} \times \mathbb{R}_+ \times \mathbb{R}_+)$ such that for every $(\xi_0, z_0) \in \mathbb{R} \times \mathbb{R}$, $d \in L^\infty(\mathbb{R}_+; \mathbb{R})$, $\theta \in L^\infty(\mathbb{R}_+; \mathbb{R})$ the unique solution of (111) with $(\xi(0), z(0)) = (\xi_0, z_0)$ that corresponds to $d \in L^\infty(\mathbb{R}_+; \mathbb{R})$, $\theta \in L^\infty(\mathbb{R}_+; \mathbb{R})$ exists for all $t \geq 0$ and satisfies the following estimates for all $t \geq 0$:

$$\limsup_{t \to +\infty}(|\xi(t)|) \leq \sqrt{2\varepsilon}, \tag{112}$$

$$|\xi(t)| \leq \exp(-ct/2)|\xi_0| + \sqrt{\frac{2a}{c(1+\exp(z_0))}} \left( \|d\|_\infty + \left( \|\theta\|_\infty - 1 - \exp(z_0) \right)^+ \right), \quad (113)$$

$$\limsup_{t \to +\infty} (|\xi(t)|)$$
$$\leq \sqrt{\frac{2a}{c\left(1+\exp\left(\lim_{t \to +\infty}(z(t))\right)\right)}} \left( \limsup_{t \to +\infty}(|d(t)|) + \left( \limsup_{t \to +\infty}(|\theta(t)|) - 1 - \exp\left(\lim_{t \to +\infty}(z(t))\right) \right)^+ \right) \quad (114)$$

as well as estimates (83), (82). When $\theta \in \mathbb{R}$ is a constant, the stability properties of system (111) with output $y = \xi$ can be stated with precision:

- estimates (83), (82) guarantee the p-UBIBS property,
- estimate (112) guarantees the zero p-OAG property with asymptotic residual constant equal to $\sqrt{2\varepsilon}$,
- estimate (113) guarantees the p-IOS property with gain function $\gamma(s) = s\sqrt{2a/c}$ and residual constant $\alpha = (|\theta| - 1)^+ \sqrt{2a/c}$.

When $\theta \in \mathbb{R}$ is a constant and in the disturbance-free case (i.e., when $d \equiv 0$), the stability properties of system (111) with output $y = \xi$ can also be stated with precision:

- estimates (83), (82) guarantee the Lagrange stability property,
- the set $\{(0, z)' \in \mathbb{R}^2, z \in \mathbb{R}\}$ is invariant for (111); however, it should be noticed that (depending on the value of $\theta \in \mathbb{R}$) there may be equilibrium points of system (111) out of the set $\{(0, z)' \in \mathbb{R}^2, z \in \mathbb{R}\}$,
- estimate (113) guarantees the p-UGAOS property and residual constant $\alpha = (|\theta| - 1)^+ \sqrt{2a/c}$,
- estimates (112), (113) show that system (111) with output $y = \xi$ is p-GOA with asymptotic residual constant equal to $\tilde{\alpha} = \min\left(\sqrt{2\varepsilon}, (|\theta| - 1)^+ \sqrt{2a/c}\right) \leq \alpha$,
- estimates (112) and (114) guarantee the estimate $\limsup_{t \to +\infty}(|\xi(t)|) \leq \min\left(\sqrt{2\varepsilon}, \sqrt{\frac{2a}{c(1+l)}}(|\theta| - 1 - l)^+\right)$, where $l = \exp\left(\lim_{t \to +\infty}(z(t))\right)$; thus, we can obtain an explicit condition for the convergence of the output to zero, namely, $\lim_{t \to +\infty}(\xi(t)) = 0$ when $|\theta| \leq 1 + l$. Consequently, we get the following result:

"At least one of the following holds:
(i) $\lim_{t \to +\infty}(\xi(t)) = 0$,
(ii) $|\theta| > 1$ and $z(t) \leq \lim_{s \to +\infty}(z(s)) < \ln(|\theta| - 1)$, for all $t \geq 0$."

It is clear that system (111) with output $y = \xi$ is rich of stability properties. ◁

We end this section by providing sufficient conditions for the p-OAG property. The following theorem is a novel result which is close in spirit to LaSalle's theorem (see [9]) which is not applicable to systems with inputs. In order to overcome the presence of inputs the following result exploits a family of functions (instead of a single function).

**Theorem 9 (LaSalle-like conditions for the p-OAG property):** *Suppose that for every $x_0 \in \mathbb{R}^n$ and $d \in L^\infty(\mathbb{R}_+; D)$ the solution of (36) with initial condition $x(0) = x_0$ is bounded. Moreover, suppose that there exist a non-decreasing function $b \in C^0(\mathbb{R}_+; \mathbb{R}_+)$ and families of functions $H_s \in C^1(\mathbb{R}^n)$, $Q_s \in C^0(\mathbb{R}^n; \mathbb{R}_+)$ being locally Lipschitz, parameterized by $s \geq 0$, such that the following properties hold for all $s \geq 0$, $x \in \mathbb{R}^n$:*

$$Q_s(x) = 0 \Rightarrow |h(x)| \leq b(s) \tag{115}$$

$$\max\{\nabla H_s(x) f(x,d) : d \in D, |d| \leq s\} \leq -Q_s(x) \tag{116}$$

*Then system (36), (2) satisfies the p-OAG property with asymptotic residual constant $\tilde{\alpha} = b(0)$ and asymptotic gain function $\tilde{\gamma}(s) = b(s) - b(0)$.*

Before we proceed with the proof of Theorem 9 it is worth noticing some things about Theorem 9. First, the family of functions $H_s \in C^1(\mathbb{R}^n)$ parameterized by $s \geq 0$ is not a family of positive definite (or even non-negative) functions. Secondly, we say that $H_s \in C^1(\mathbb{R}^n)$, $Q_s \in C^0(\mathbb{R}^n; \mathbb{R}_+)$ are families of functions $H_s \in C^1(\mathbb{R}^n)$, $Q_s \in C^0(\mathbb{R}^n; \mathbb{R}_+)$, parameterized by $s \geq 0$, instead of saying that $H, Q$ are functions of $(x,s) \in \mathbb{R}^n \times \mathbb{R}_+$ because there is no regularity requirement with respect to $s \geq 0$; the functions $H_s \in C^1(\mathbb{R}^n)$, $Q_s \in C^0(\mathbb{R}^n; \mathbb{R}_+)$ are not assumed to be continuous (or even measurable) with respect to $s \geq 0$. When (60) and

$$V(x) \geq \varphi(|d|) \Rightarrow \dot{V}(x,d) = \nabla V(x) f(x,d) \leq -c(x,|d|) \tag{117}$$

hold for certain functions $a \in K_\infty$, $V \in C^1(\mathbb{R}^n; \mathbb{R}_+)$, a non-decreasing, continuous function $\varphi : \mathbb{R}_+ \to \mathbb{R}_+$ and a locally Lipschitz function $c : \mathbb{R}^n \times \mathbb{R}_+ \to \mathbb{R}_+$ with $c(x,s) > 0$ for all $(x,s) \in \mathbb{R}^n \times \mathbb{R}_+$ with $V(x) > \varphi(s)$, then (115), (116) hold with

$$H_s(x) = \frac{1}{2}\left((V(x) - \varphi(s))^+\right)^2, \; Q_s(x) = (V(x) - \varphi(s))^+ \min\{c(x,l) : 0 \leq l \leq s\} \text{ and } b(s) = a^{-1}(\varphi(s))$$

**Proof of Theorem 9:** Let arbitrary $x_0 \in \mathbb{R}^n$, $d \in L^\infty(\mathbb{R}_+; D)$ be given and consider the bounded solution $x(t) \in \mathbb{R}^n$ of (36) with initial condition $x(0) = x_0$. Define $s = \|d\|_\infty$. Since the solution is bounded, there exists $R \in \mathbb{R}$ such that $H_s(x(t)) \geq R$ for all $t \geq 0$. Exploiting (116) and the fact that $s = \|d\|_\infty$ we get for all $t \geq 0$:

$$\int_0^t Q_s(x(\tau))d\tau \leq H_s(x(0)) - H_s(x(t)) \leq H_s(x(0)) - R \tag{118}$$

Moreover, the facts that $d \in L^\infty(\mathbb{R}_+; D)$, $x(t) \in \mathbb{R}^n$ is bounded and (36) imply that $\dot{x}(t)$ is of class $L^\infty(\mathbb{R}_+)$. Since $Q_s \in C^0(\mathbb{R}^n; \mathbb{R}_+)$ is locally Lipschitz, it follows that the function $t \mapsto Q_s(x(t))$ is uniformly continuous. From Barbălat's Lemma (see [5, 9]), we have:

$$\lim_{t \to +\infty}(Q_s(x(t))) = 0 \tag{119}$$

We next claim that

$$\limsup_{t \to +\infty}(|y(t)|) \leq b(\|d\|_\infty) \tag{120}$$

We prove (120) by contradiction. Suppose that $\limsup_{t \to +\infty}(|y(t)|) > b(\|d\|_\infty)$. Then by virtue of (2) there exists $\varepsilon > 0$ and a sequence of times $\{t_i \geq 0 : i = 1, 2, ...\}$ with $\lim_{i \to +\infty}(t_i) = +\infty$ such that $|y(t_i)| = |h(x(t_i))| \geq b(\|d\|_\infty) + \varepsilon$. Since the solution $x(t) \in \mathbb{R}^n$ is bounded, it follows that the sequence $\{x(t_i) \in \mathbb{R}^n : i = 1, 2, ...\}$ is bounded. By taking a subsequence if necessary, we may assume that there exists $x^* \in \mathbb{R}^n$ such that $\lim_{i \to +\infty}(x(t_i)) = x^*$. Consequently, by exploiting continuity of $h$, $Q_s$, (119) and the fact that $s = \|d\|_\infty$ we get $|h(x^*)| \geq b(\|d\|_\infty) + \varepsilon = b(s) + \varepsilon > b(s)$ and $Q_s(x^*) = 0$. Clearly, this contradicts (115).

Inequality (120) and the fact that $x_0 \in \mathbb{R}^n$, $d \in L^\infty(\mathbb{R}_+; D)$ are arbitrary imply the p-OAG property with asymptotic residual constant $\tilde{\alpha} = b(0)$ and asymptotic gain function $\tilde{\gamma}(s) = b(s) - b(0)$. The proof is complete. ◁

**Example 10:** Consider the following planar system

$$\begin{aligned} \dot{\xi} &= 0 \\ \dot{y} &= -yp(\xi, y, d) + (1 - y^2)^+ q(\xi, y, d) \\ x &= (\xi, y)' \in \mathbb{R}^2, d \in \mathbb{R} \end{aligned} \tag{121}$$

where $p, q : \mathbb{R}^3 \to \mathbb{R}$ are locally Lipschitz functions that satisfy $q(0,0,0) = 0$ and $p(\xi, y, d) > 0$ for all $(\xi, y, d) \in \mathbb{R}^3$. System (121) is not necessarily a system that satisfies the IOS property (e.g., consider the input-free case $p(\xi, y, d) \equiv 1$ and $q(\xi, y, d) = 2y$). However, system (121) satisfies the p-UBIBS property (use Theorem 5 with $H_1(x) = \xi^2/2$, $H_2(x) = y^2/2$ and $R(s) \equiv 1/2$).

The function
$$V(\xi, y) = y^2/2 \qquad (122)$$

gives for all $(x,d) \in \mathbb{R}^2 \times \mathbb{R}$:

$$\dot{V}(\xi, y, d) = -y^2 p(\xi, y, d) + y(1-y^2)^+ q(\xi, y, d)$$

It follows that (117) holds with $\varphi(s) = 1/2$ and $c(\xi, y, |d|) = y^2 (V(\xi, y) - 1/2)^+ \min(p(\xi, y, |d|), p(\xi, y, -|d|))$. Moreover, (60) holds with $a(s) = s^2/2$. Consequently, (115), (116) hold with

$$H_s(x) = \left((V(x) - \varphi(s))^+\right)^2 / 2, \quad Q_s(x) = y^2 (V(\xi, y) - 1/2)^+ \min\{p(\xi, y, l): |l| \leq s\} \text{ and } b(s) \equiv 1$$

We conclude from Theorem 9 that the OAG property

$$\limsup_{t \to +\infty} (|y(t)|) \leq 1$$

holds for all $d \in L^\infty(\mathbb{R}_+; \mathbb{R})$. ◁

**Example 11:** Consider again system (111), which was studied in Example 9 and was shown (by means of Theorem 8) to be a system that satisfies the p-UBIBS property. Here we show the validity of inequality (112) by means of Theorem 9. We consider system (111) with output $y = \xi$ and we apply Theorem 9 with

$$H_s(x) = -\exp(z), \quad Q_s(x) = \Gamma\left(\frac{1}{2}\xi^2 - \varepsilon\right)^+ \text{ and } b(s) \equiv \sqrt{2\varepsilon}$$

Inequality (112), i.e., the zero p-OAG property with asymptotic residual constant equal to $\sqrt{2\varepsilon}$, is a direct consequence of Theorem 9.

Notice that in this case the function $H_s(x)$ has negative values. ◁

## 7) The Need of Output Stability Properties for Systems with Disturbances in Adaptive Control

Since disturbances are always present in every system and since disturbances appear as external inputs for a system in the standard mathematical formalism that is used in control theory, it becomes clear that there are no control systems without external inputs. When adaptive controllers are used, one has to study the stability properties of the closed-loop system with external inputs. This study is not easy. In this section we focus on the simplest system, namely the scalar system

$$\dot{y} = \theta y + u + d$$
$$y \in \mathbb{R}, \theta \in \mathbb{R}, u \in \mathbb{R}, d \in \mathbb{R} \tag{124}$$

where $\theta \in \mathbb{R}$ is an unknown parameter of the system, $u \in \mathbb{R}$ is the control input and $d \in \mathbb{R}$ is the (unknown and possibly time-varying) external disturbance.

A nonlinear adaptive controller may be designed based on the methodologies proposed in [10]. Indeed, the adaptive controller

$$u = -cy - \hat{\theta} y - q y^3$$
$$\frac{d\hat{\theta}}{dt} = \Gamma y^2 \tag{125}$$

where $c, \Gamma > 0$, $q \geq 0$ are the controller parameters (gains) guarantees a number of stability properties in the disturbance-free case (i.e., when $d \equiv 0$). Indeed, by defining $z = \hat{\theta} - \theta$, we obtain the equivalent version of the closed-loop system (124) with (125):

$$\dot{y} = -cy - zy - qy^3 + d$$
$$\dot{z} = \Gamma y^2 \tag{126}$$
$$x = (y, z)' \in \mathbb{R}^2, d \in \mathbb{R}$$

The radially unbounded function $V(x) = \frac{1}{2} y^2 + \frac{1}{2\Gamma} z^2$ satisfies the following inequality for all $y, z, d \in \mathbb{R}$:

$$\dot{V}(x, d) = -cy^2 - qy^4 + yd \leq -\frac{c}{2} y^2 - qy^4 + \frac{1}{2c} d^2 \tag{127}$$

Consequently, the disturbance-free version of system (126):

- by virtue of Theorem 1 and Proposition 1, satisfies the Lagrange and Lyapunov stability properties,
- the set $\{(0, z)' \in \mathbb{R} \times \mathbb{R}, z \in \mathbb{R}\}$ is invariant for (126); it is a set of equilibrium points for system (126),
- by virtue of Theorem 2 is GAOS.

However, we cannot show UGAOS for system (126). Still, the adaptive controller (125) works nicely when disturbances are absent. Figure 1 shows the plots of the state variables and the control input in the disturbance-free case.

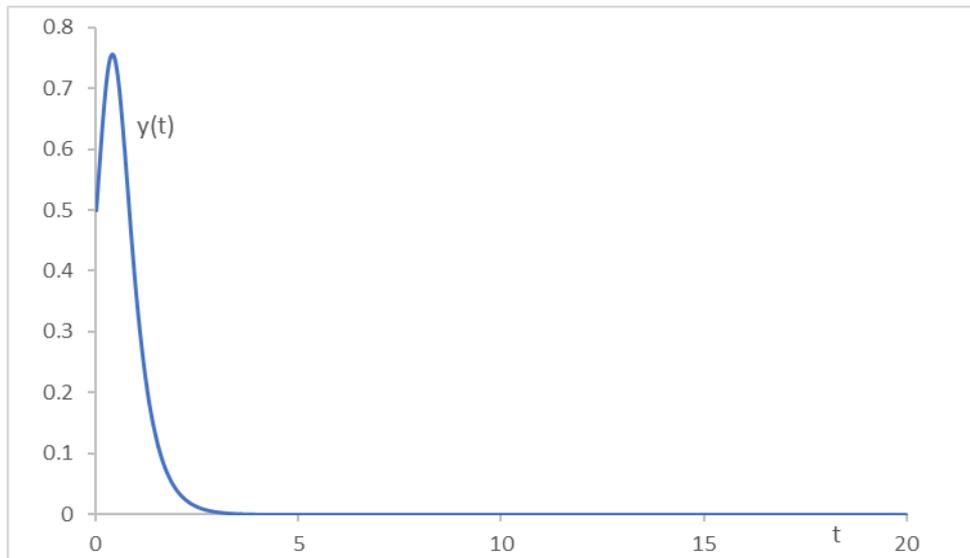

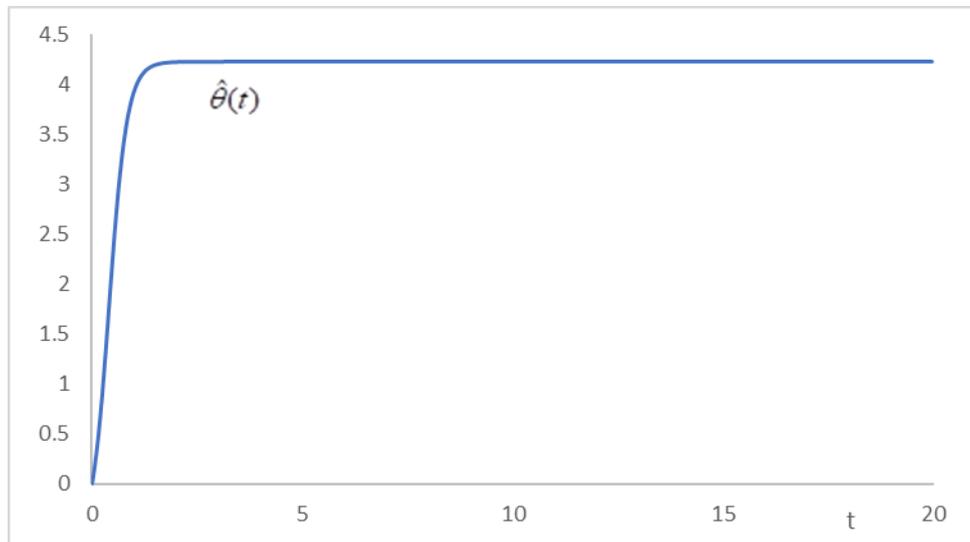

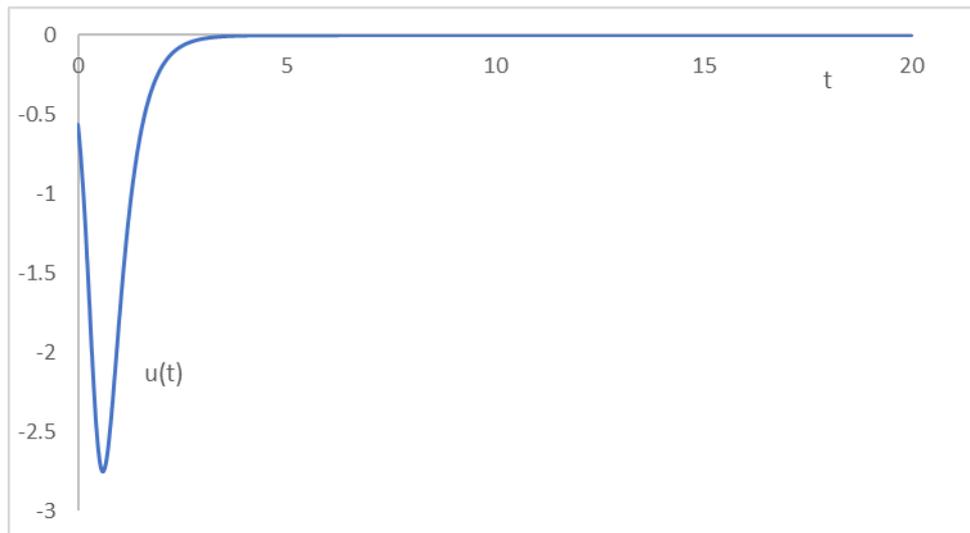

**Fig. 1:** The evolution of the state variables and the control input for the closed-loop system (124), (125), with $d \equiv 0$, $\theta = 3$, $c = 1, \Gamma = 10, q = 0.5$ and $y(0) = 0.5$, $\hat{\theta}(0) = 0$.

Another advantage of the adaptive controller (125) is its simplicity. The adaptive controller (125) consists of two parts: the static part $u = -cy - \hat{\theta}y - qy^3$ and the dynamic part $\frac{d\hat{\theta}}{dt} = \Gamma y^2$. The dynamic part of the adaptive controller (125) is called "an identifier" that produces the parameter estimate $\hat{\theta}(t)$ of the true value of the parameter $\theta$.

By virtue of Theorem 4 and inequality (127) we can also guarantee that system (126) is forward complete. However, we cannot guarantee boundedness of the state variables for arbitrary inputs $d \in L^\infty(\mathbb{R}_+;\mathbb{R})$. In order to guarantee state boundedness through estimate (127), one has to assume that the disturbance has only "finite energy", i.e., that $d \in L^2(\mathbb{R}_+;\mathbb{R}) \cap L^\infty(\mathbb{R}_+;\mathbb{R})$. The reason that we cannot prove boundedness of the state variables for arbitrary inputs $d \in L^\infty(\mathbb{R}_+;\mathbb{R})$ is that the states are not bounded. Figure 2 shows the problem in the case of a persistent disturbance $d \in L^\infty(\mathbb{R}_+;\mathbb{R})$ with $d \notin L^2(\mathbb{R}_+;\mathbb{R})$.

Figure 1 and Figure 2 show another aspect of adaptive control schemes which are similar to (125). The identifier fails to guarantee convergence of the parameter estimate $\hat{\theta}(t)$ to $\theta$. The variable $\hat{\theta}$ works more or less as a dynamic controller gain which can only increase with respect to time. The increase of $\hat{\theta}(t)$ does not stop in the case of a persistent disturbance $d \in L^\infty(\mathbb{R}_+;\mathbb{R})$ with $d \notin L^2(\mathbb{R}_+;\mathbb{R})$ and this causes a slow increase of the magnitude of the control input $|u(t)|$.

The problems shown in Figure 2 have been known in the control community for a long time now and many modifications have been proposed for handling such situations (see [2, 3]). The most popular modification is the so-called "σ-modification" (see [18] for the application of σ-modification even to infinite-dimensional systems) where only the identifier is modified by adding an extra term of the form "$-\sigma\hat{\theta}$" where $\sigma > 0$ is a constant; this term tries to keep $\hat{\theta}(t)$ bounded. Therefore, the adaptive controller (125) with σ-modification reads

$$\begin{aligned} u &= -cy - \hat{\theta}y - qy^3 \\ \frac{d\hat{\theta}}{dt} &= \Gamma y^2 - \sigma\hat{\theta} \end{aligned} \quad (128)$$

By defining $z = \hat{\theta} - \theta$, we obtain the equivalent version of the closed-loop system (124) with (128):

$$\begin{aligned} \dot{y} &= -cy - zy - qy^3 + d \\ \dot{z} &= \Gamma y^2 - \sigma z - \sigma\theta \\ x &= (y,z)' \in \mathbb{R}^2, d \in \mathbb{R} \end{aligned} \quad (129)$$

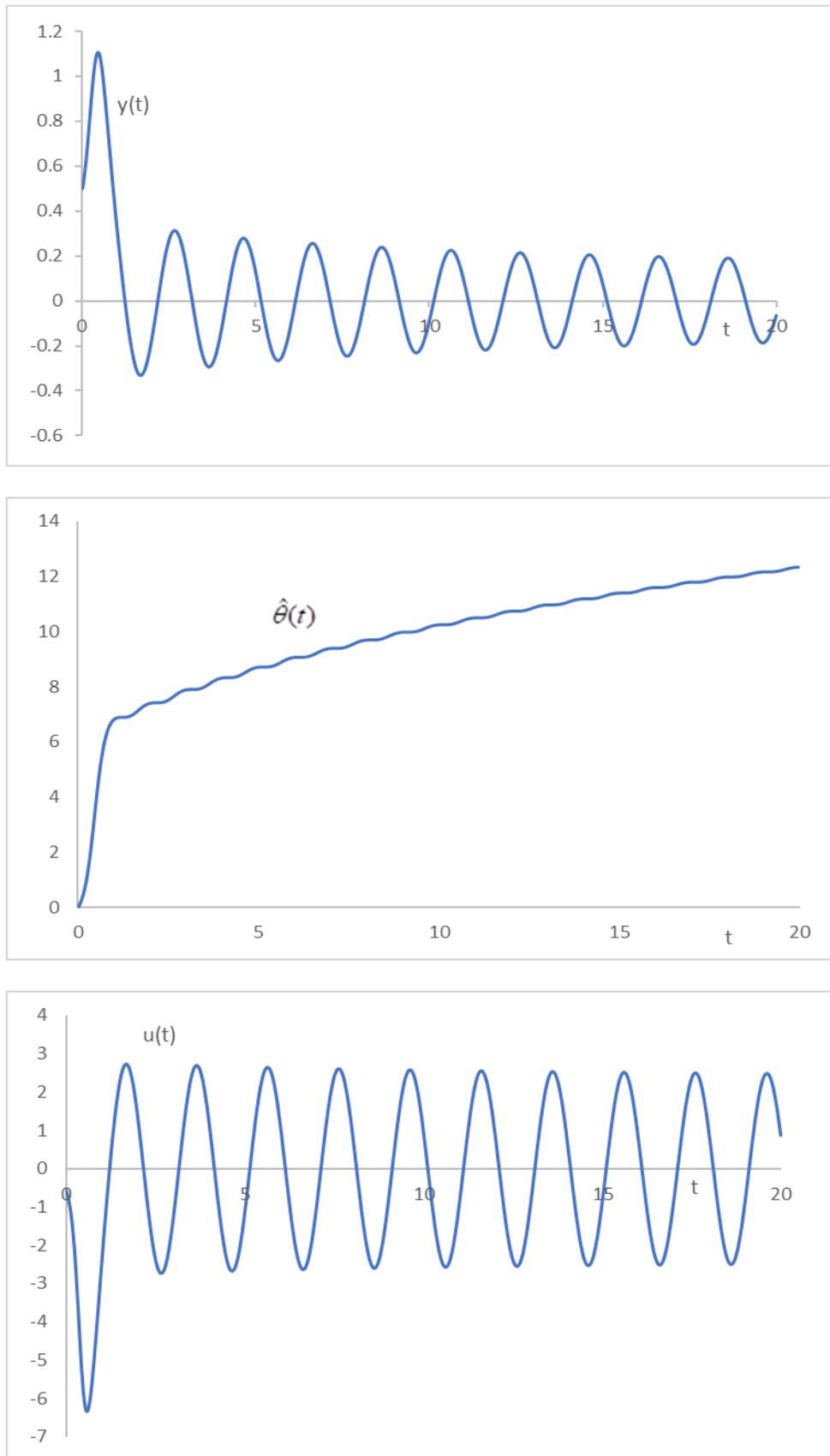

**Fig. 2:** The evolution of the state variables and the control input for the closed-loop system (124), (125), with $d(t) = 2\sin(\pi t)$, $\theta = 3$, $c = 1, \Gamma = 10, q = 0.5$ and $y(0) = 0.5$, $\hat{\theta}(0) = 0$.

Again, an advantage of the adaptive controller (128) is its simplicity. The adaptive controller (128) consists again of two parts: the static part $u = -cy - \hat{\theta}y - qy^3$ and the dynamic part $\dfrac{d\hat{\theta}}{dt} = \Gamma y^2 - \sigma \hat{\theta}$. Due to the σ-modification we cannot expect the dynamic part of the adaptive controller (128) to work as "an identifier" that produces a reliable parameter estimate $\hat{\theta}(t)$ of the true value of the parameter $\theta$. Instead, as noted above, the variable $\hat{\theta}$ works more or less as a dynamic controller gain which cannot increase too much.

The radially unbounded function $V(x) = \dfrac{1}{2}y^2 + \dfrac{1}{2\Gamma}z^2$ for system (129) satisfies the following inequality for all $y, z, d \in \mathbb{R}$:

$$\dot{V}(x,d) = -cy^2 - qy^4 - \dfrac{\sigma}{\Gamma}z^2 - \dfrac{\sigma}{\Gamma}z\theta + yd$$

$$\leq -\dfrac{c}{2}y^2 - qy^4 - \dfrac{\sigma}{2\Gamma}z^2 + \dfrac{\sigma}{2\Gamma}\theta^2 + \dfrac{1}{2c}d^2 \tag{130}$$

$$\leq -\min(c,\sigma)V(x) + \dfrac{\sigma}{2\Gamma}\theta^2 + \dfrac{1}{2c}d^2$$

Integrating the above differential inequality, we obtain the following estimate for every $d \in L^\infty(\mathbb{R}_+;\mathbb{R})$:

$$|x(t)| \leq \exp(-\min(c,\sigma)t/2)\sqrt{\max(\Gamma,1/\Gamma)}|x(0)|$$
$$+ \sqrt{\dfrac{\max(\Gamma,1)}{c\min(c,\sigma)}}\|d\|_\infty + |\theta|\sqrt{\dfrac{\sigma}{\min(\Gamma,1)\min(c,\sigma)}} \tag{131}$$

Estimate (131) guarantees the p-ISS property for system (129) with gain function $\gamma(s) = s\sqrt{\dfrac{\max(\Gamma,1)}{c\min(c,\sigma)}}$ and residual constant $\alpha = |\theta|\sqrt{\dfrac{\sigma}{\min(\Gamma,1)\min(c,\sigma)}}$. In the disturbance-free case we guarantee p-UGAS with residual constant $\alpha = |\theta|\sqrt{\dfrac{\sigma}{\min(\Gamma,1)\min(c,\sigma)}}$.

Although the p-ISS property is a very strong stability property and estimate (131) also guarantees exponential convergence, there is a problem here. The problem is that the residual constant is proportional to the absolute value of the unknown parameter, i.e., $\alpha \sim |\theta|$. One may say that the stability analysis that we have performed is very conservative and that is why we get such a result. Unfortunately, this is not the case. If we perform an analysis of the equilibrium points of system (129) then we can easily see that there are equilibrium points of the form $x = (y,z) = \left(\pm\sqrt{\dfrac{\sigma(\theta-c)}{\Gamma+\sigma q}}, -\dfrac{\Gamma c + \theta\sigma q}{\Gamma+\sigma q}\right)$ when $\theta > c$. Consequently, the residual constant (which must be large enough so that the ball with radius the residual constant contains all equilibrium points) must be approximately proportional to the absolute value of the unknown parameter. The problem is real and it is shown in Figure 3 for the disturbance-free case and in Figure 4 for the case of a persistent disturbance.

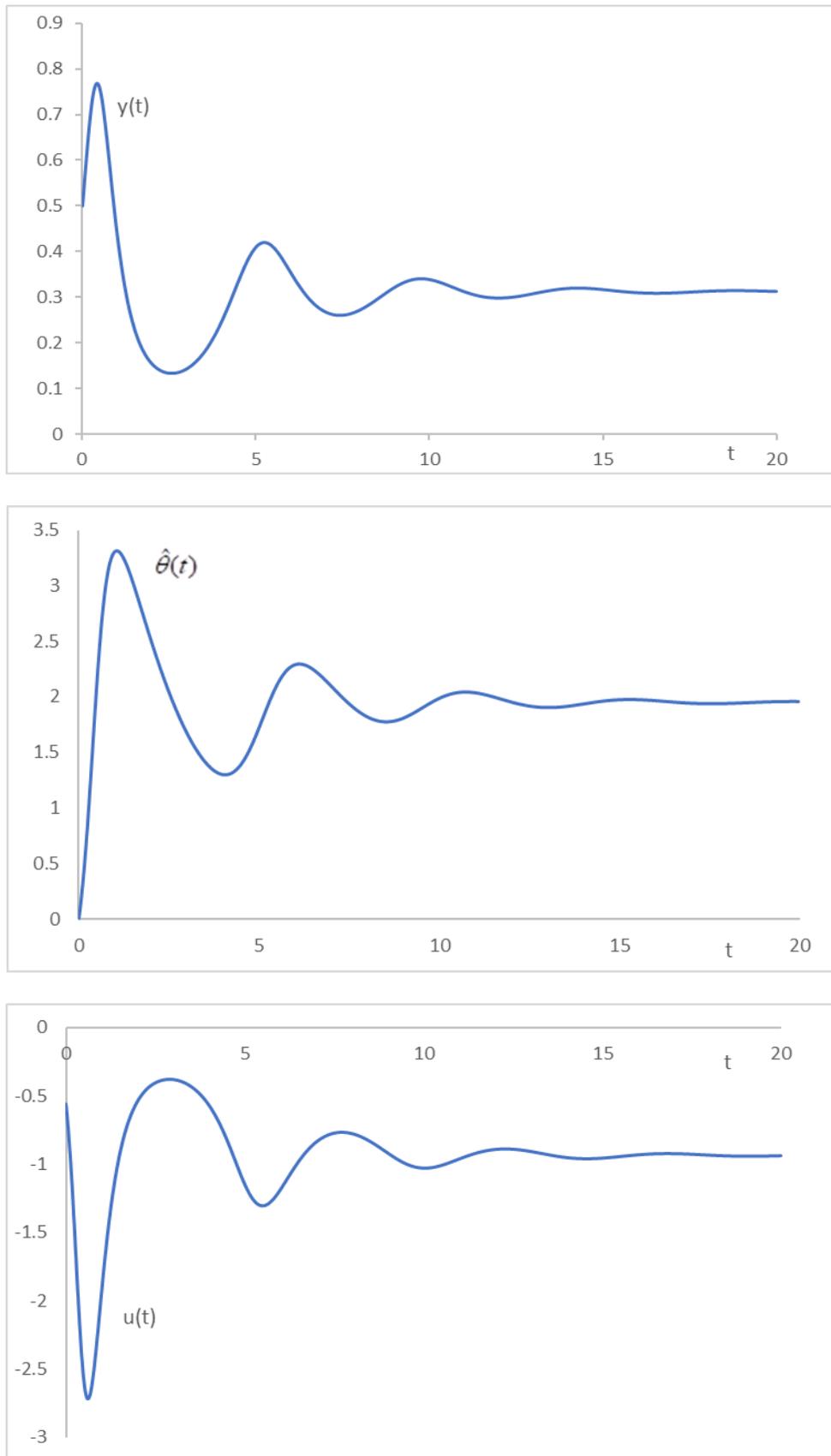

**Fig. 3:** The evolution of the state variables and the control input for the closed-loop system (124), (128), with $d(t) \equiv 0$, $\theta = 3$, $c = 1, \Gamma = 10, \sigma = q = 0.5$ and $y(0) = 0.5$, $\hat{\theta}(0) = 0$.

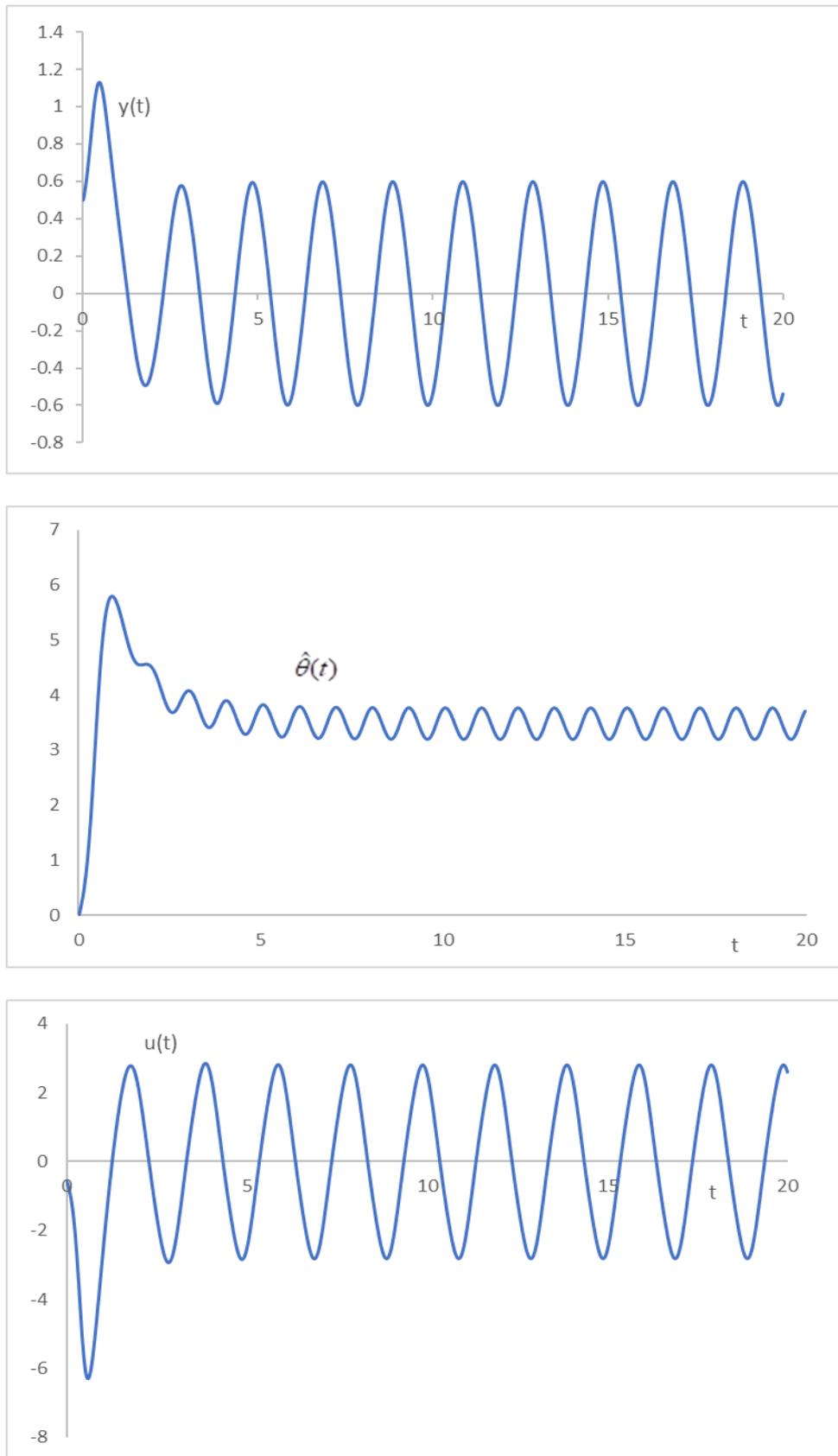

**Fig. 4:** The evolution of the state variables and the control input for the closed-loop system (124), (128), with $d(t) = 2\sin(\pi t)$, $\theta = 3$, $c = 1, \Gamma = 10, \sigma = q = 0.5$ and $y(0) = 0.5$, $\hat{\theta}(0) = 0$.

The advantages of σ-modification are its simplicity, the fact that we can guarantee exponential convergence in the disturbance-free case with assignable convergence rate, the fact that we can guarantee the p-ISS property with assignable gain function. The disadvantage of σ-modification is the fact that the residual constant is not a-priori known due to the dependence of the residual constant on the unknown parameter $\theta$. Therefore, we cannot control the offset in the response of the closed-loop system -even in the disturbance-free case- and that can give a poor performance.

Having identified the weaknesses of the σ-modification (poor transients, regulation to an unknown neighborhood of the origin dependent on the unknown parameter), we remind the reader of Example 9, which introduced an adaptive system that overcomes these deficiencies of the σ-modification. This example, with its contrasting tools relative to the σ-modification (deadzone, update rate that decreases as the adaptive gain grows, several forms of "nonlinear damping") is generalized from the scalar linear case in [6, 8].

We end the section by noticing the precision that the stability properties (and the problems) of the closed-loop systems under adaptive controllers are described by using the output stability notions that we introduced in the present set of notes.

**Notes and comments**

- There are many other stability properties beyond the ones that are discussed in the present work. We have focused mainly on the stability properties that are useful for adaptive control. We hope the reader appreciates the clear presentation of this collection of stability properties, with a scope slightly expanded beyond adaptive control.

- Output stability properties are very much related to the so-called "partial stability" properties. See for instance [27]. Output stability is a special case of stability with respect to two measures; see [11, 12, 16, 26].

- Controller (125) with $q=0$ is a "high-gain" adaptive controller (see [14, 15, 28]).

- ISS was introduced in control theory by E. D. Sontag in [19]. IOS was introduced in control theory by E. D. Sontag and Y. Wang in [23]. Both notions were studied by E. D. Sontag and Y. Wang; see [19, 20, 21, 22, 23, 24, 25]. The role that IOS and ISS played in modern feedback control theory is huge. The UBIBS property, i.e., property (40) with $\bar{\alpha}=0$ was introduced in [24].

- Some of the results contained in the present set of notes are novel results. For example, Theorem 4, Theorem 5, Theorem 8 and Theorem 9 have not appeared elsewhere. A less general version of Theorem 8 first appeared in [8].

- While writing this set of notes, Prof. Antoine Chaillet mentioned to the authors that the same proof procedure that is used in the proof of Theorem 1.2 in [24] can be exploited in order to obtain the following sufficient condition for p-IOS:

**Theorem 10 (Lyapunov conditions for p-IOS):** *Suppose that (36) is p-UBIBS and there exist a function $a \in K_\infty$, a non-decreasing, continuous function $\varphi: \mathbb{R}_+ \to \mathbb{R}_+$, a function $V \in C^1(\mathbb{R}^n; \mathbb{R}_+)$ with $V(0)=0$ and a function $\sigma \in KL$ for which the following inequality holds for all $x \in \mathbb{R}^n$*

$$a(|h(x)|) \leq V(x) \tag{132}$$

*and the following implication holds for all $x \in \mathbb{R}^n$, $d \in D$:*

$$V(x) \geq \varphi(|d|) \Rightarrow \dot{V}(x,d) = \nabla V(x) f(x,d) \leq -\sigma(V(x), |x|) \tag{133}$$

*Then system (36), (2) satisfies the p-IOS property. Moreover, if $\varphi(0)=0$ then system (36), (2) satisfies the IOS property.*

In other words, it is possible to replace the UBIBS property in Theorem 1.2 in [24] by the less demanding p-UBIBS property.

- It is well known that, for time-invariant systems described by ODEs with locally Lipschitz right-hand sides, GAS and UGAS are actually equivalent properties; see [29]. However, GAOS is not equivalent to UGAOS for arbitrary output maps.

**Acknowledgments:** The authors would like to Professor Antoine Chaillet for all the discussions in the early stages of writing this document.